\documentclass [12pt]{amsart}

\usepackage{amssymb,amsxtra,amsfonts}
\usepackage{epsfig}             %
\usepackage{graphics}

\usepackage{verbatim}  
\usepackage{bbding}   

\newenvironment{ppb}[1]
{\ \!\!\!\!\!\!\!\!\!\!\!\!\!\!\!\!\!\!\!\!\!\!\!\!\!\!\!\!\!\!\!\!\!\!\!\!\!\!\!\! {\bf PPB------------------------------------------------------------------------------------------------PPB}\newline \tiny {#1}
\  \newline\normalsize\phantom{f}\!\!\!\!\!\!\!\!\!\!\!\!\!\!\!\!\!\!\!\!\!\!\!\!\!\!\!\!\!\!\!\!\!\!\!\!\!\!\!\! {\bf PPB------------------------------------------------------------------------------------------------PPB}\newline}{}

\def\reE@DeclareMathSymbol#1#2#3#4{%
    \let#1=\undefined
    \DeclareMathSymbol{#1}{#2}{#3}{#4}}
\DeclareSymbolFont{symbolsC}{U}{txsyc}{m}{n}
\SetSymbolFont{symbolsC}{bold}{U}{txsyc}{bx}{n}
\DeclareFontSubstitution{U}{txsyc}{m}{n}
\reE@DeclareMathSymbol{\strictiff}{\mathrel}{symbolsC}{76}

\newcommand\beq{\begin{equation}}
\newcommand\eeq{\end{equation}}
\newcommand\bal{\begin{align*}}
\newcommand\eal{\end{align*}}   
\newcommand\bmx{\left(\begin{matrix}}
\newcommand\emx{\end{matrix}\right)}
\newcommand\bsmx{\left(\begin{smallmatrix}}
\newcommand\esmx{\end{smallmatrix}\right)}

\newcommand{\st}{\ \bigl\vert\ }

\def\part#1{\frac{\partial\phantom{q}}{\partial#1}}

\newcommand {\flb}{\lbrack\!\lbrack}
\newcommand {\frb}{\rbrack\!\rbrack}
\newcommand {\flp}{(\!(}
\newcommand {\frp}{)\!)}

\newcommand{\glue}[1]{\underset{#1}{\strictiff}}
\newcommand{\fus}{\circledast}

\newcommand{\MDR}{\mathcal{M}_{\text{\rm DR}}}
\newcommand{\MB}{\mathcal{M}_{\text{\rm B}}}

\newcommand{\MDol}{\mathcal{M}_{\text{\rm Dol}}}

\newcommand{\Jac}{\text{\rm Jac}}

\newcommand{\HH}{\text{\rm H}}
\newcommand{\hh}{\mathop{\rm H}}

\newcommand{\Lie}{{\mathop{\rm Lie}}}

\DeclareMathOperator{\ISto}{{\IS}to} 

\DeclareMathOperator{\Sect}{\mathop{\rm Sect}}

\newcommand{\Sym}{\mathop{\rm Sym}}
\newcommand{\Fun}{\mathop{\rm Fun}}

\newcommand{\papk}[3]{  \ _{#1}^{\phantom{#3}}\cA_{#2}^{#3}     }
\newcommand{\gahr}{\papk{G}{H}{r}}

\newcommand{\Ad}{{\mathop{\rm Ad}}}
\newcommand{\ad}{{\mathop{\rm ad}}}

\DeclareMathOperator{\pr}{pr}

\newcommand{\Prod}{\prod}

\DeclareMathOperator{\Hom}{Hom}         

\newcommand{\SL}{{\mathop{\rm SL}}}

\newcommand{\GL}{{\mathop{\rm GL}}}

\newcommand{\U}{{\rm {U}}}	

\newcommand{\irr}{{\rm irr}}

\newcommand{\End}{\mathop{\rm End}}

\newcommand{\reg}{{\mathop{\footnotesize\rm reg}}}

\newcommand{\hk}{{hyperk\"ahler }}   

\newcommand{\PVI}{{$\text{\rm P}_{\text{\rm VI}}$}}   

\newcommand{\bH}{{\bf H}}

\newcommand{\IA}{\mathbb{A}}
\newcommand{\IB}{\mathbb{B}}
\newcommand{\IC}{\mathbb{C}}

\newcommand{\IM}{\mathbb{M}}

\newcommand{\IP}{\mathbb{P}}                                     
\newcommand{\IQ}{\mathbb{Q}}                           
\newcommand{\IR}{\mathbb{R}}                           
\newcommand{\IS}{\mathbb{S}}

\newcommand{\IZ}{\mathbb{Z}}

\newcommand{\cA}{\mathcal{A}}
\newcommand{\cB}{\mathcal{B}}
\newcommand{\cC}{\mathcal{C}}

\newcommand{\cD}{\mathcal{D}}

\newcommand{\cG}{\mathcal{G}}

\newcommand{\cM}{\mathcal{M}}

\newcommand{\cO}{\mathcal{O}}

\newcommand{\cR}{\mathcal{R}}

\newcommand{\cU}{\mathcal{U}}

\newcommand{\gM}{       \mathfrak{M}     }

\newcommand{\g}{       \mathfrak{g}     }

\newcommand{\lt}{\mathfrak{t}}

\newcommand{\gl}{       \mathfrak{gl}     } 

\newcommand{\lp}{\mathfrak{p}}

\newcommand{\lu}{\mathfrak{u}}

\newcommand{\wt}{\widetilde}

\newcommand{\wh}{\widehat}

\newcommand{\al}{\alpha}

\newcommand{\ga}{\gamma}
\newcommand{\de}{\delta}

\newcommand{\Ga}{\Gamma}

\newcommand{\La}{\Lambda}

\newcommand{\si}{\sigma}

\newcommand{\Si}{\Sigma}
\renewcommand{\th}{\theta}

\makeatletter
 \newlength{\typesize}
\makeatother

\newlength{\vvoff}
\newlength{\hhoff}

\newcommand{\locateoffcenter}[1]{%
\addtolength{\vvoff}{-0.25\typesize}%
\raisebox{\vvoff}{\hspace{\hhoff}\makebox(0,0){\smash{#1}}}
}
\newcommand{\object}[1]{%
\setlength{\vvoff}{0pt}%
\setlength{\hhoff}{0pt}%
\locateoffcenter{#1}
}

\newcommand{\swlabel}[1]{%
\setlength{\vvoff}{-0.5\typesize}%
\setlength{\hhoff}{0.75\typesize}%
\locateoffcenter{#1}
}
\newcommand{\nwlabel}[1]{%
\setlength{\vvoff}{-0.5\typesize}%
\setlength{\hhoff}{-0.75\typesize}%
\locateoffcenter{#1}
}
\newcommand{\selabel}[1]{%
\setlength{\vvoff}{0.5\typesize}%
\setlength{\hhoff}{0.75\typesize}%
\locateoffcenter{#1}
}

\def\mapright#1{\smash{
        \mathop{\longrightarrow}\limits^{#1}}}

\def\mapdown#1{\Big\downarrow
        \rlap{$\vcenter{\hbox{$\scriptstyle#1$}}$}}

\def\underset#1#2{\ \smash{\mathop{ #2 }\limits_{#1}}\ }

\newcommand{\pf}{\begin{bpf}}

\newcommand{\pfms}{\begin{bpfms}}
\newcommand{\epf}{\end{bpf}\hfill$\square$\\}           
\newcommand{\epfms}{\end{bpfms}\hfill$\square$\\}       

\newcommand{\idea}{\begin{bidea}}

\newcommand{\eidea}{\end{bidea}\hfill$\square$\\}           

\newcommand{\sk}{\begin{bsk}}    

\newcommand{\esk}{\end{bsk}\hfill$\square$\\}           
\newcommand{\sketch}{\begin{bsketch}}

\newcommand{\esketch}{\end{bsketch}\hfill$\square$\\}

\newtheorem {hypo}{\bf\hspace{-\parindent}Hypothesis}
\newtheorem {thm}[hypo]{Theorem}   

\newtheorem {cor}[hypo]{Corollary}

\theoremstyle{remark}

\renewenvironment{ppb}[1]{}{}

\begin{document}

\thispagestyle{empty}

\begin{center}
{\Large \bf TH\`ESE D'HABILITATION DE L'UNIVERSIT\'E PARIS XI}

\ 

\ 

Sp\'ecialit\'e : MATH\'EMATIQUES

\ 

\ 

pr\'esent\'ee par

\ 

\ 

{\bf \Large Philip BOALCH}

\end{center}

\ 

\ 

\

\noindent
Sujet de la th\`ese :

\ 

\begin{center}

{\bf \large
Geometry of moduli spaces of meromorphic connections on curves, Stokes data, wild nonabelian Hodge theory, hyperk\"ahler manifolds, isomonodromic deformations, Painlev\'e equations,  and relations to Lie theory.
}
\end{center}

\vskip 4cm

\begin{center}
Soutenue le 12/12/12 devant le jury compos\'e de
\end{center}

\ppb{$$
\begin{alignat*}{3}
\text{1} && \text{2}\\
\text{1} && \text{2}
\end{alignat*}
$$}
 $$\begin{array}{lcl}
\text{ALEKSEEV, Anton} &\quad& \text{rapporteur externe}\\
\text{BOST, Jean-Beno$\hat{\text{\i}}$t} && \text{}\\

\text{HITCHIN, Nigel} && \text{}\\
\text{SABBAH, Claude} && \text{rapporteur interne}\\
\text{SCHIFFMANN, Olivier} && \text{}\\
\text{SIMPSON, Carlos} && \text{rapporteur externe}\\
\end{array}
$$

\newpage
\thispagestyle{empty}

\ 

\vskip 5cm

{\em Adresse de l'auteur} : 

\ 

Philip BOALCH 

DMA, \'Ecole Normale Sup\'erieure, 

45 rue d'Ulm, 

75005 Paris, France

\ 

{boalch@dma.ens.fr}

www.math.ens.fr/$\sim$boalch

\newpage
\thispagestyle{empty}

\begin{center}
{\Large \bf Table des mati\`eres}
\end{center}

\ 

\noindent
{\large \bf I.  Pr\'esentation des travaux}

1 Introduction\hfill \pageref{sn: intro}

2 Hyperk\"ahler moduli spaces and wild nonabelian Hodge theory \hfill \pageref{sn: wnabh}

3 The nonlinear Schwarz's list \hfill \pageref{sn: nlsl}

4 Dual exponential maps and the geometry of quantum groups \hfill \pageref{sn: dexp}

5 Fission \hfill\pageref{sn: fission}

6 Wild character varieties \hfill \pageref{sn: wcvs}

7 Braid group actions from isomonodromy \hfill\pageref{sn: bafi}

8 ``Logahoric'' connections on parahoric bundles\hfill\pageref{sn: logahoric}

9 Dynkin diagrams for isomonodromy systems \hfill\pageref{sn: slims}

\

\noindent
{\large \bf II.  Articles}

\noindent
\emph{Stokes matrices, Poisson Lie groups and Frobenius manifolds}

\hspace{1em}Invent. Math. \textbf{146} (2001) 479--506 

\noindent
\emph{G-bundles, isomonodromy and quantum {W}eyl groups} 

\hspace{1em}Int. Math.  Res. Not. (2002), no.~22, 1129--1166

\noindent
\emph{Wild non-abelian {H}odge theory on curves} (avec O. Biquard)
 
\hspace{1em}Compositio Math. \textbf{140}, no.~1 (2004) 179--204

\noindent
\emph{Painlev\'e equations and complex reflections}

\hspace{1em}Ann. Inst. Fourier \textbf{53}, no.~4 (2003) 1009--1022

\noindent
\emph{From {K}lein to {P}ainlev\'e via {F}ourier, {L}aplace and {J}imbo}

\hspace{1em}Proc. London Math. Soc. \textbf{90}, no.~3 (2005), 167--208

\noindent
\emph{The fifty-two icosahedral solutions to {P}ainlev\'e {VI}}

\hspace{1em}J. Reine Angew. Math. \textbf{596} (2006) 183--214

\noindent
\emph{Some explicit solutions to the {R}iemann--{H}ilbert
  problem}

\hspace{1em}%
IRMA Lect. Math. Theor. Phys. \textbf{9} (2006) 85--112

\noindent
\emph{{Higher genus icosahedral {P}ainlev\'e curves}}

\hspace{1em}Funk. Ekvac. (Kobe) \textbf{50} (2007) 19--32

\noindent
\emph{Quasi-{H}amiltonian geometry of meromorphic connections}

\hspace{1em}Duke Math. J. \textbf{139}, no.~2 (2007) 369--405

\noindent
\emph{Regge and Okamoto symmetries}

\hspace{1em}Comm. Math. Phys. \textbf{276} (2007) 117--130

\noindent
\emph{Quivers and difference Painlev\'e equations}

\hspace{1em}%
CRM Proc. Lecture Notes, \textbf{47} (2009) 25--51

\noindent
\emph{Through the analytic halo: Fission via irregular singularities}

\hspace{1em}Ann. Inst. Fourier  \textbf{59}, no.~7 (2009) 2669--2684

\noindent
\emph{Riemann--Hilbert for tame complex parahoric connections} 

\hspace{1em}Transform. groups \textbf{16}, no.~1 (2011) 27--50

\noindent
\emph{Simply-laced isomonodromy systems}

\hspace{1em}Publ. Math. I.H.E.S. \textbf{116}, no.~1 (2012) 1-68 %

\noindent
\emph{Geometry and braiding of Stokes data; Fission and wild character varieties}

\hspace{1em}Annals of Math., to appear (accepted 6/11/12)

 \ 

\ 

\noindent
{\large Articles originaux non-pres\'ent\'es}

\ 

\noindent
P.~P. Boalch, \emph{{S}ymplectic manifolds and isomonodromic deformations}, Adv. in  Math. \textbf{163} (2001), 137--205.

\noindent
\bysame, \emph{Irregular connections and Kac-Moody root systems},  arXiv:0806.1050, June 2008, 31pp.  (largely subsumed in the last two articles presented)

\  

\noindent
{\large Articles de vulgarisation}

\

\noindent
[A] P.~P. Boalch, \emph{Brief introduction to Painlev\'e VI}, SMF, S\'eminaires et congr\`es, vol 13 (2006) 69-78.

\noindent
[B] \bysame, \emph{Six results on Painlev\'e VI}, 
SMF, S\'eminaires et congr\`es, vol 14, (2006) 1-20. %

\noindent
[C] \bysame, \emph{Towards a nonlinear Schwarz's list},
In: The many facets of geometry: a tribute to Nigel Hitchin.
J-P. Bourguignon, O. Garcia-Prada and S. Salamon (eds), OUP (2010) pp. 210-236, (arXiv:0707.3375 July 2007).

\noindent
[D] \bysame, \emph{Noncompact complex symplectic and hyperk\"ahler manifolds},
Notes for M2 cours specialis\'e 2009, 75pp., www.math.ens.fr/$\sim$boalch/hk.html

\noindent
[E] \bysame, \emph{Hyperk\"ahler manifolds and nonabelian Hodge theory on (irregular) curves}, 16pp. 2012, arXiv:1203.6607

\newpage

\noindent {\bf\large 1. Introduction.} \label{sn: intro}

The aim of this manuscript is to outline the main work the author has done since 1999.

The principal theme is the study of the geometry of certain moduli spaces attached to smooth complex algebraic curves, and the nonlinear differential equations that naturally arise when the curve, and some other parameters, are varied.

For example given a curve $\Si$ one may consider the Jacobian variety $\Jac(\Si)$ which may be viewed as the moduli space of degree zero holomorphic line bundles on $\Si$.
This gives a map 
$$\Si\mapsto \Jac(\Si)$$
associating an abelian variety to an algebraic curve.
If we now vary the curve in a family we obtain a family of abelian varieties.

In work of Weil, Mumford, Narasimhan--Seshadri  and others it was understood that there is a similar picture for higher rank vector bundles, in effect replacing the structure group $\IC^*$ appearing in the line bundle case by the non-abelian group $\GL_n(\IC)$, provided one introduces a stability condition (which is automatic in the line bundle case) 
\cite{Weil38, mumford.icm62, NarSes}.
This gives a map 
$$\Si\mapsto \cU_n(\Si)$$
where $\cU_n(\Si)$ denotes the moduli space of stable degree zero holomorphic vector bundles on $\Si$, a non-abelian analogue of the Jacobian $\Jac(\Si)\cong\cU_1(\Si)$. The theorem of Narasimhan--Seshadri says that $\cU_n(\Si)$ is homeomorphic to  $\Hom^\irr(\pi_1(\Si),U_n)/U_n$, the space of irreducible unitary  representations of the fundamental group of $\Si$.

In work of Hitchin, Simpson and others it was understood that for many purposes it is better to consider a ``complexified version'' of this story.
They defined the notion of Higgs bundle, which consists of a vector bundle 
$E\to\Si$ together with a Higgs field 
$\Phi\in \Gamma(\Omega^1 \otimes\End E)$. 
The moduli space $\MDol(\Si,n)$ of stable rank $n$ degree zero Higgs bundles is then a partial compactification of the cotangent bundle $T^*\cU_n(\Si)$ of the space of stable bundles and there is a diffeomorphism
$$\MDol(\Si,n)\cong \MDR(\Si,n)$$
with the moduli space $\MDR(\Si,n)$ of (stable) holomorphic connections on rank $n$ vector bundles $V\to \Si$.
This isomorphism may be interpreted both as a nonabelian analogue of Hodge theory  \cite{Sim-hbls} (noting that the nonabelian cohomology space $\HH^1(\Si,\GL_n(\IC))$ classifies rank $n$ vector bundles with flat connection),
and as a rotation of complex structure on an underlying  hyperk\"ahler manifold %
\cite{Hit-sde}. 
In turn the Riemann--Hilbert correspondence, taking a flat connection to its monodromy representation, yields an analytic isomorphism
 $$\MDR(\Si,n)\  \cong\  \MB(\Si,n) := \Hom^\irr(\pi_1(\Si),\GL_n(\IC))/\GL_n(\IC)$$
to the space of irreducible complex representations of the fundamental group of $\Si$.

If we now vary the curve $\Si$ in a family over a base $\IB$ then the De\! Rham and Betti spaces $\MDR(\Si,n), \MB(\Si,n)$ fit together into fibre bundles over $\IB$, both of which admit natural flat (Ehresmann/nonlinear) connections on their total space (and correspond to each other via the Riemann--Hilbert correspondence).
This nonlinear connection is the nonabelian analogue of the Gauss--Manin connection (see Simpson \cite{Sim94ab}).
The parallel  between this  and earlier work on isomonodromic deformations was pointed out in \cite{smid}.

\begin{figure}[h] 	
\centerline{\input{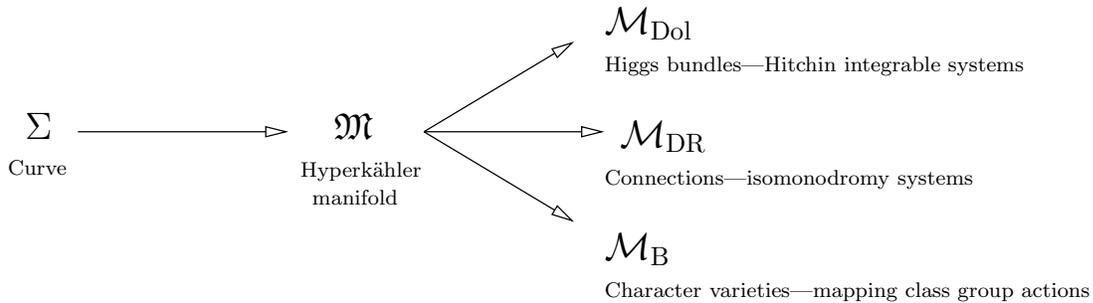}}	
\caption{Basic setup (from the survey \cite{ihptalk}).} \label{fig: basic setup}
\end{figure}

Most of the work to be presented in this manuscript is concerned with various aspects of the extension of this story
when one considers meromorphic connections, rather than just the holomorphic connections appearing above (in the definition of $\MDR$), and especially (but not exclusively) the case of meromorphic connections with irregular singularities.

On one hand this yields many more moduli spaces, and even the case when the underlying curve is the Riemann sphere is now extremely interesting.
For example new complete \hk four-manifolds (gravitational instantons) arise as the simplest nontrivial examples of such moduli spaces (see \S2 below). 

Secondly in this meromorphic case the nonlinear connections sometimes have explicit descriptions (and a longer history, cf. \cite{schles-icm1908}) and in many examples the resulting nonlinear differential equations have actually appeared in both physical and mathematical problems.
A basic set of examples are the Painlev\'e equations, which appear when the moduli spaces have complex dimension two. They are experiencing something of a renaissance since their appearance in high energy physics and string theory. 

Thirdly there are extra deformation parameters which occur in the case of irregular meromorphic connections, beyond the moduli of the underlying curve with marked points. These extra parameters (controlling the ``irregular type'' of the connection) nonetheless behave exactly like the moduli of the curve and similarly lead to nonlinear braid group actions on the moduli spaces. 
(Interestingly if one considers meromorphic connections on $G$-bundles this brings the $G$-braid groups into play \cite{bafi}.)

Fourthly the Betti description of irregular connections is more complicated than  the fundamental group representation appearing above; it  involves  ``Stokes data'' enriching the fundamental group representation. One theme of this work has been to understand the geometry of such spaces of  Stokes data and their relation to other parts of mathematics.

\newpage
\noindent%
{\bf\large2.\,Hyperk\"ahler moduli spaces \& wild nonabelian Hodge theory} \label{sn: wnabh}

The article \cite{wnabh} with O. Biquard extends the nonabelian Hodge correspondence, between Higgs bundles and local systems (or flat connections), to a correspondence between meromorphic Higgs bundles and (irregular) meromorphic connections on smooth complex algebraic curves.
In particular this constructs a large class of complete \hk manifolds.
This work is surveyed in \cite{ihptalk}.

\noindent{\bf Background. }
In \cite{Hit-sde, Don87} Hitchin and Donaldson established a correspondence between stable Higgs bundles and local systems (or holomorphic connections) on a smooth compact complex algebraic curve.
This was extended (to higher dimensional projective varieties, and higher rank structure groups) by Corlette and Simpson \cite{Cor88, Sim-hbls}, who also interpreted this correspondence as a nonabelian analogue of Hodge theory (cf. \cite{Simpson-nabh, simp-hfnc}). 
In Hitchin's framework the correspondence arises naturally from \hk geometry: one simply rotates to a different complex structure in the \hk family to move from the Higgs bundle moduli space to the moduli space of connections.
The \hk viewpoint was extended to the moduli spaces which arise in the case of higher dimensional projective varieties by Fujiki \cite{fujiki-hk}, but he also noted (\cite{fujiki-hk} p.3) that in fact all the moduli spaces which arise in this way embed into a moduli space that arises in the case of a curve\footnote{More pointedly (and recently), Simpson \cite{Sim04} p.2  stated: ``the
irreducible components of moduli varieties of flat connexions which are
known, are all isomorphic to moduli varieties of representations on curves''.}.
To get {\em new} moduli spaces one may consider meromorphic connections. 
The nonabelian Hodge correspondence was extended by Simpson \cite{Sim-hboncc} to the case of meromorphic connections on open curves satisfying a tameness assumption (so the resulting moduli spaces are basically representations of the fundamental group of the curve).
On the other hand meromorphic Higgs bundles on curves 
(with arbitrary poles) 
had been considered algebraically  
 \cite{ AHH, Bea, Nit-higgs, Bot, Mar} and it was shown that 
they had many of the properties of the nonsingular case, such as being fibred by abelian varieties/admitting the structure of algebraically completely integrable Hamiltonian system (the meromorphic Hitchin integrable systems).

\noindent{\bf Main result.}
The main result of \cite{wnabh}, which is reviewed succinctly in \cite{ihptalk}, 
can be summarised as follows.
Fix a general linear group $\GL_n(\IC)$ and consider a smooth compact algebraic curve $\Si$ with some marked points $a_1,\ldots a_m\in \Si$.
At each point choose an `irregular type' $Q_i$, some weights $\th_i$ and  a residue element $\tau_i+\si_i+N_i$, in the notation of \cite{ihptalk}.
This data determines a moduli space $\MDR(\Si,\th,\tau,\si, N)$
of isomorphism classes of stable meromorphic connections with compatible parabolic structures and the given irregular types, weights and residue orbits.
Similarly one may choose data $Q',\th',\tau',\si', N'$ and consider a moduli space $\MDol(\Si,\th',\tau',\si', N')$
of stable meromorphic Higgs bundles with compatible parabolic structures and the given irregular types, weights and residue orbits.

\begin{thm}[\cite{wnabh}]\label{thm: main}
The moduli space $\MDR(\Si,\th,\tau,\si, N)$ 
of meromorphic connections is a \hk manifold and it is
 naturally diffeomorphic to the moduli space
\newline
\noindent
$\MDol(\Si,\th',\tau',\si', N')$
of meromorphic Higgs bundles if the data are related as follows:
$$Q_i'=-Q_i/2,\quad N'_i=N_i,\quad \th_i' = -\tau_i-[-\tau_i],\quad \tau_i' = -(\tau_i+\th_i)/2,\quad\si_i' = -\si_i/2$$

\noindent
where $[\,\cdot\,]$ denotes the (component-wise) integer part.
Moreover the \hk metrics are {\em complete} if the nilpotent parts are zero ($N=0$) and there are no strictly semistable objects, and this may be ensured by taking the parameters to be off of some explicit hyperplanes (\cite{wnabh} \S8.1).
\end{thm}

This correspondence is established by passing through solutions to Hitchin's self-duality equations, and the map from meromorphic connections to such solutions 
(i.e. constructing a harmonic metric for irregular connections) was established earlier by Sabbah \cite{Sab99} in the case of trivial Betti weights (this is the irregular analogue of the result of Donaldson and Corlette).
The approach of \cite{wnabh} is simpler, due to a `straightening trick' avoiding the Stokes phenomenon and enabling a simpler construction of initial metric, leading to the full correspondence and the construction of the moduli spaces.
(The \hk quotient of \cite{wnabh} is a strengthening of the complex symplectic quotient description of $\MDR$ in the irregular case \cite{Boa, smid}, which  used a similar straightening trick.) 

It should be emphasised perhaps that this construction gives new examples of complete \hk manifolds even in complex dimension two (i.e. real four manifolds), cf. \cite{ihptalk} \S3.2---these are referred to as ``gravitational instantons'' by physicists, and Atiyah \cite{atiyah-hk} has emphasised their purely mathematical significance, as the quaternionic analogue of algebraic curves.

\noindent
{\bf Further developments.}
Witten
\cite{witten-wild} has used these hyperk\"ahler manifolds to extend his approach to the geometric Langlands correspondence to the wildly ramified case, extending his work with Kapustin \cite{kapustin-witten} 
and with Gukov \cite{gukov-witten}.
Other physicists have been very interested in trying to construct such \hk 
metrics in a more explicit fashion and have related the existence of such metrics to the so-called Kontsevich--Soibelman wall-crossing formula (see e.g. \cite{gmn09}).
(As far as I know there is still no rigorous example of 
such an explicit approach.)
See for example \cite{wildquivergt, xie12} for more on the role in high energy physics of these \hk moduli spaces of irregular singular solutions to Hitchin's equations.
Within mathematics, T. Mochizuki has extended some aspects of the wild nonabelian Hodge correspondence to higher dimensions, and has used this to prove a conjecture of Kashiwara \cite{mochizuki-wild}.

\

\newpage

\noindent%
{\large\bf3. The nonlinear Schwarz's list} \label{sn: nlsl}

In  the articles \cite{pecr, k2p, icosa, octa, ipc} the author discovered,  classified and constructed many algebraic solutions of the 
sixth Painlev\'e equation.
This work is surveyed in \cite{nlsl}.

\noindent
{\bf Background.}
The classical list of Schwarz \cite{Schwarz} is a list of the algebraic solutions of the Gauss hypergeometric equation, and they are related to the finite subgroups of $\SL_2(\IC)$.
This Gauss hypergeometric equation is a linear differential equation and it is the simplest explicit example of a Gauss--Manin connection.
The {\em nonabelian} Gauss--Manin connections are natural 
nonlinear connections which arise when one considers the nonabelian cohomology of a family of varieties.
The simplest example is the family of Painlev\'e VI differential equations, which arises from $\HH^1(X,G)$ with $X$ a four-punctured Riemann sphere, and $G=\SL_2(\IC)$.
The ``nonlinear'' analogue of  Schwarz's list is thus a list of algebraic solutions of the sixth Painlev\'e equation.
The nonlinear case is  considerably more difficult 
since: 
a) there is no simple a priori finiteness: for example it is not enough to go through finite subgroups of $\SL_2(\IC)$, 
b) even if a solution is proven to exist, it is still highly nontrivial to actually construct it (one needs to explicitly solve a family of {\em nonrigid} Riemann--Hilbert problems),
c) the affine Weyl group of type $F_4$ acts on the set of algebraic solutions, so one needs to be careful that any ``new'' solution is not just a transformation of a known solution. %

One motivation is that nonlinear differential equations such as 
Painlev\'e VI arise in many nonlinear problems in geometry and high energy physics, and it is known that most solutions of Painlev\'e VI are new transcendental functions, not expressible in terms of simpler special functions. 
One often finds very special geometric objects correspond to the special explicit algebraic solutions.
Thus for example Hitchin \cite{Hit-tei} constructs some four-dimensional Einstein manifolds from some special algebraic solutions, and Dubrovin \cite{Dub95long-with.ApE.note} Appendix E, relates certain algebraic solutions of Painlev\'e VI to certain algebraic Frobenius manifolds ($=$ mathematical TQFTs).

\noindent{\bf Previous results.}
Before working on this project there were explicit algebraic solutions constructed by Hitchin \cite{Hit-Poncelet,Hit-tei, Hit-Ed},
Dubrovin \cite{Dub95long-with.ApE.note} Appendix E, Dubrovin--Mazzocco \cite{DubMaz00}, 
and Kitaev/Andreev were writing a series of papers (\cite{Kit-sfit6,And-Kit-CMP}) containing many algebraic solutions.

\noindent{\bf Main results.}
In brief there are three continuous families (due to Okamoto, Hitchin and Dubrovin), one discrete family (due to Picard and Hitchin) of algebraic solutions to Painlev\'e VI, and
then:

\begin{thm}[\cite{pecr, k2p, icosa, octa, ipc, nlsl}] 
There are at least $45$ inequivalent exceptional/sporadic algebraic solutions of the Painlev\'e VI differential equation.
\end{thm}

Nine of these sporadic solutions are not due to the author\footnote{Such counting is difficult due to the $W_{\text{aff}}(F_4)$ action: $1$ solution is due to Andreev--Kitaev, $1$ to Dubrovin, $2$ to Dubrovin--Mazzocco, $5$ to Kitaev, cf. 
\cite{P6survey} p.18, \cite{nlsl}. 
This is corroborated in \cite{lis-tyk}.}. 
The other $36$ solutions were found and constructed explicitly in \cite{pecr, k2p, icosa, octa, ipc} (the number $45$ appears in the last section of \cite{nlsl}---see also \cite{P6survey}).

Particular highlights of this result include:

1) There is a solution \cite{pecr, k2p} coming from Klein's simple group of order $168$, and it is not related to any finite subgroup of $\SL_2(\IC)$,

2) Most of the solutions are related to the symmetry groups of the platonic solids and all such ``platonic'' solutions were classified in \cite{icosa, octa} (and the outstanding platonic solutions were constructed in these articles and \cite{ipc}).
Note that $19$ of the $52$ icosahedral solutions are not sporadic, as explained in \cite{icosa}.

3) There is an icosahedral solution (\cite{icosa} Theorem B) which is ``generic'' in the sense that its parameters lie of {\em none} of the affine $F_4$ reflection hyperplanes,

4) There is a uniquely determined algebraic curve of genus $7$ canonically attached to the icosahedron (on which the largest, degree $72$, icosahedral solution is defined). An explicit plane model for this curve is as follows (\cite{ipc}):

$$
9\,(p^6\,q^2+p^2\,q^6)+
18\,p^4\,q^4+$$$$
4\,(p^6+q^6)+
26\,(p^4\,q^2+p^2\,q^4)+
8\,(p^4+q^4)+
57\,p^2\,q^2+$$$$
20\,(p^2+q^2)+
16=0.
$$
$$\text{The genus seven icosahedral Painlev\'e curve.}$$

5) The degree $18$ solution of Dubrovin--Mazzocco, which involved an elliptic curve that took many pages of $40$ digit integers to write down (in the preprint version of \cite{DubMaz00} on the arXiv), has a simple parameterisation 
(\cite{icosa} Theorem C), and the underlying elliptic curve may be given by the formula $u^2=s(8s^2-11s+8)$.

\noindent{\bf Further developments.}
By 2005-6 there seemed to be nowhere left to look for more algebraic solutions, so the list of known solutions was lectured about in 2006 \cite{P6survey}, written up in \cite{nlsl}, and the problem of proving there were no more solutions was set (last page of \cite{hit60} or \cite{P6survey}).
In 2008 Lisovyy and Tykhyy \cite{lis-tyk} showed (by computer calculation) 
that there are no more algebraic solutions---in particular the count of $45$ exceptional solutions is as in \cite{nlsl}.
(Their article is also highly recommended for the wonderful colour pictures illustrating  the topological structure of the branching of the solutions.)

{\tiny
\noindent{\bf Other perspectives.}
As part of this work on algebraic solutions of the nonlinear Painlev\'e VI equation, many algebraic solutions of certain {(nonrigid)} linear differential equations were found.
This problem is of interest in its own right (see e.g. Baldassarri--Dwork \cite{Bald-Dwork-algsols} %
and the literature on the Grothendieck--Katz conjecture, such as Katz \cite{Katz-Algsols}).
For example the $52$ icosahedral solutions of Painlev\'e VI in 
\cite{icosa} constitute an extension of the icosahedral part of Schwarz's classical list to the case of rank two Fuchsian systems with four poles on the Riemann sphere---the first $10$ rows correspond to the $10$ icosahedral rows on Schwarz's list. (In principal these are all pullbacks of hypergeometric equations, but in practice the pullbacks are hard to compute a priori: this approach was pursued in \cite{Chuck1} and \cite{And-Kit-CMP}, but few new inequivalent solutions were constructed). See also \cite{octa} for the octahedral and tetrahedral cases, and e.g. \cite{srops} \S3 for some explicit rank three connections with finite monodromy generated by three reflections. 
Except for \cite{B-vdW}, previous extensions of Schwarz's list, such as \cite{BH89}, remained in the world of rigid differential equations---things are then much easier as there are no ``accessory parameters'' (the moduli spaces are zero dimensional).
}

\newpage
\noindent%
{\bf \large4.\,Dual exponential maps and the geometry of quantum groups} \label{sn: dexp}

The articles \cite{smapg, bafi} defined and studied a natural class of holomorphic maps on the dual of the Lie algebra of any complex 
reductive group. 
This arose from a moduli-theoretic realisation of the classical limit of the Drinfeld--Jimbo quantum group---in other words the author discovered that their quantum group {\em is} a quantisation of a very simple moduli space of irregular connections.
Corollaries include  new direct proofs of theorems of Kostant, Duistermaat and Ginzburg--Weinstein and a geometric understanding of the so-called quantum Weyl group.
(These results were proved for $G=\GL_n(\IC)$ in \cite{smapg} and extended to other complex reductive groups in \cite{bafi}---here we mainly restrict to $\GL_n(\IC)$ for simplicity.)

\noindent{\bf Background.}
Let $G=\GL_n(\IC)$ and let $\g = \End(\IC^n)$ denote its Lie algebra, so that the dual vector space $\g^*$ is naturally a complex Poisson manifold. Using the trace pairing 
the space $\g^*$ is identified with $\g$, so that $\g$ inherits a complex Poisson structure and the symplectic leaf through $A\in \g$ is its adjoint orbit. 

\begin{figure}[h] 	
\centerline{\input{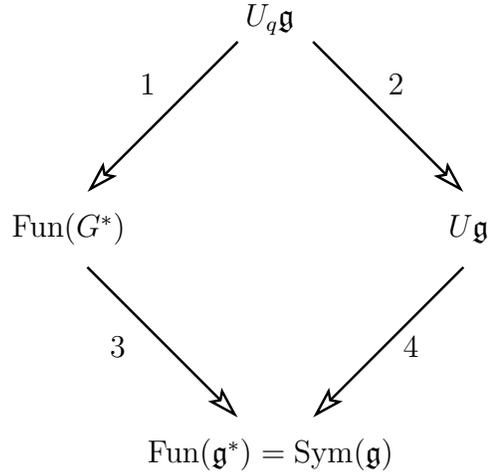}}	
\caption{Simplifying the Drinfeld--Jimbo quantum group.} \label{fig: hopf algebras}
\end{figure}

In the theory of quantum groups the main Poisson manifolds which appear are certain nonlinear analogues of the linear Poisson manifolds $\g^*$. 
The most important example is the following (it is due to Drinfeld/Semenov-Tian-Shansky
\cite{Drin86} example 3.2 in infinitesimal form, \cite{STS83},\cite{DKP} p.185, \cite{AlekMalk94} p.170).
Let $B_+\subset G$ be a Borel subgroup (such as the upper triangular matrices) and let $T\subset B_+$ be a maximal torus (such as the diagonal matrices). Let $B_-\subset G$ be the opposite Borel (so that $B_-\cap B_+ = T$), and let 
$\delta : B_\pm\to T$
be the natural projection (taking the ``diagonal part'').
The {\em standard dual Poisson Lie group} of $G$ is 
$$G^* = \{ (b_-,b_+)\in B_-\times B_+\st \de(b_-)\de(b_+) = 1\}\subset G\times G$$
which is an algebraic group of the same dimension as $G$.
Sometimes it is convenient to consider the universal cover of $G^*$, by including an element $\La\in\lt=\Lie(T)$ such that 
$\de(b_\pm)=\exp(\pm\pi i\La)$, although the resulting group is no longer algebraic.
The group $G^*$ 
admits a natural Poisson structure, which may be defined geometrically (see \cite{smapg} \S2, following \cite{LuRatiu}).
The symplectic leaves of $G^*$ are obtained by fixing the 
conjugacy class of the product
\beq\label{eq: sp leaf condn}
b_-^{-1}b_+\in G.
\eeq
The relevance (and importance) of the Poisson manifold $G^*$ is that the Drinfeld--Jimbo quantum group is a deformation quantisation of it. 
More precisely, there is the following diagram
of Hopf algebras.
To understand this first recall that the algebra of functions on a Lie group is a commutative Hopf algebra, which is cocommutative if and only if the underlying Lie group is abelian;  a ``quantum group'' is a non-commutative  non-cocommutative Hopf algebra. Thus one may simplify a quantum group in various ways (see Figure \ref{fig: hopf algebras}).

\ppb{
\begin{equation*}
 \setlength{\unitlength}{36pt}
 \begin{picture}(2.2,2.2)(0,0)
 \put(2,1){\object{$U\g$}}
 \put(0,1){\object{$\Fun(G^*)$}}
 \put(1,2){\object{$U_q\g$}}
 \put(1,0){\object{$\Fun(\g^*)=\Sym(\g)$}}
 \put(0.75,1.75){\vector(-1,-1){0.5}}
 \put(0.6,1.8){\nwlabel{$1$}}
 \put(1.25,1.75){\vector(1,-1){0.5}}
 \put(1.4,1.8){\swlabel{$2$}}
 \put(1.75,0.75){\vector(-1,-1){0.5}}
 \put(1.4,0.2){\selabel{$4$}}
 \put(0.25,0.75){\vector(1,-1){0.5}}
 \put(0.1,.5){\swlabel{$3$}}
 \end{picture}
\end{equation*}
}

\newpage 
Here $U_q\g$ is the Drinfeld--Jimbo quantum group (which is a noncommutative, non-cocommutative Hopf algebra), 
$\Fun(\g^*)= \Sym(\g)$ is the Poisson algebra of functions on $\g^*$
(which is a commutative, cocommutative Hopf algebra, using the additive group structure on $\g^*$), 
$\Fun(G^*)$ is the Poisson algebra of functions on $G^*$
(which is the algebra of functions on a noncommutative group, so is a commutative, non-cocommutative Hopf algebra), 
and $U\g$ is the universal enveloping algebra of $\g$
(a non-commutative, cocommutative Hopf algebra).

Arrow 1 is due to De\,Concini--Kac--Procesi \cite{DKP},\cite{DP1565}  Theorem p.86 \S12.1:
there is an integral form of $U_q\g$ (i.e. a $\IC[q,q^{-1}]$ subalgebra) in which we can set $q=1$ and the resulting Poisson algebra is the algebra of  functions on $G^*$ (an earlier version of this result at the level of formal groups is due to Drinfeld, cf. \cite{Drin86} \S3).

Arrow 2 is the viewpoint mainly taken by Drinfeld and Jimbo 
(see \cite{Drin86} Example 6.2). 

Arrow 3 corresponds to taking the linearisation of the Poisson structure on $G^*$ at the identity, and arrow 4 corresponds to the Poincar\'e--Birkhoff--Witt isomorphism (enabling $U\g$ to be viewed as a quantisation of $\Sym\g$).

\noindent
{\bf Main Results.}
From a geometrical perspective it thus seems important to understand the Poisson manifold $G^*$.
The condition of fixing the conjugacy class of the product 
\eqref{eq: sp leaf condn}  
in order to fix a symplectic leaf is reminiscent of the condition to fix a symplectic leaf for moduli spaces of flat connections on Riemann surfaces with boundary, that one should fix the conjugacy class of the monodromy around each boundary component in order to fix a  symplectic leaf (see e.g. \cite{Aud95long}).
In fact this is not a coincidence, since:

\begin{thm}{\cite{smapg}.}
The Poisson manifold $G^*$ is isomorphic to a moduli space of meromorphic connections on the unit disk, with its natural Poisson structure, and the product \eqref{eq: sp leaf condn} is conjugate to the monodromy around the boundary circle of the connection corresponding to $(b_-,b_+)\in G^*$.
\end{thm}

Since all bundles over the disk are trivial one can write down such connections explicitly: we consider connections of the form
$$\left(\frac{A_0}{z^2} + \frac{B}{z} + holomorphic\right) dz$$
where $A_0\in \lt_\reg$ has distinct eigenvalues, and we include a framing so that in effect we only quotient by gauge transformations $g(z)$ with $g(0)=1$.
The fact that the resulting moduli space is isomorphic to $G^*$ {\em as a space} follows almost immediately from previous work on the irregular Riemann--Hilbert problem (see \cite{smapg} Theorem 5): in essence the elements of $G^*$ are the Stokes data of such connections (and the diagonal element $\La$ is the so-called ``exponent of formal monodromy'', which is just the diagonal part of $B$ in the present situation). 
The natural Poisson structure on the moduli space  is that coming from the extension of the Atiyah--Bott construction to connections with irregular singularities, from \cite{smid}.
In the present example this Poisson structure on the moduli space may be characterised quite concretely, as described in the following section.

\ 

\noindent{\bf The dual exponential map.\ }\label{ssn: dexp}
The above moduli space may be approximated by considering { global} connections on the trivial holomorphic bundle on the Riemann sphere $\IP^1(\IC)$ which have a first order pole at $\infty$ and  have the above form at $0$.
Such global connections may be written in the form 
\beq\label{eq: gconn}
\left(\frac{A_0}{z^2} + \frac{B}{z}\right) dz\eeq
for elements $B\in \g$.
If we identify $\g\cong \g^*$ using an invariant inner product then $\g$ inherits a linear complex Poisson structure from that on $\g^*$.
Thus the act of restricting such a global connection to the unit disc and taking its Stokes data yields a holomorphic map

$$\boxed{\phantom{\Bigl(}\ \ \ \nu_{A_0}: \g^*\to G^*\ \ \ \ }$$
for each choice of $A_0\in \lt_\reg$, taking an an element $B\in \g\cong\g^*$ to the Stokes data (and formal monodromy) of the corresponding connection \eqref{eq: gconn}. 
This is a highly transcendental holomorphic map between manifolds of the same dimension, and one may prove  
(cf. \cite{smapg} Lemma 31) 
it is generically a local analytic isomorphism (in particular in a neighbourhood of $0\in \g^*$).
The main result of \cite{smapg} is:

\begin{thm}
The dual exponential map $\nu_{A_0}$ is a  Poisson map for any choice of $A_0\in \lt_\reg$, relating the linear Poisson structure on $\g^*$ and the non-linear Poisson structure on $G^*$.
\end{thm}

Since equipping a vector space with a Lie bracket is equivalent to equipping the dual vector space with a linear Poisson structure, this Poisson property is evidence that $\nu_{A_0}$ should indeed be viewed as a 
dual analogue of the exponential map $\g\to G$.

In the present situation of $G=\GL_n(\IC)$ these dual exponential  maps may be directly related, using the Fourier--Laplace transform, to the Riemann-Hilbert map taking the monodromy representation of a Fuchsian system with $n+1$ poles on the Riemann sphere (see \cite{BJL81} and the exposition in \cite{k2p} \S3, \cite{nlsl} diagram 1).
This indicates how transcendental the maps $\nu_{A_0}$  are (and in particular that they are more complicated than the usual exponential map for a Lie group).
The proof given in \cite{smapg} does not use this however, and so extends to any complex reductive group, once we {\em define} Stokes data for connections on $G$ bundles (this is done in \cite{bafi}). 

Note that Drinfeld was motivated by Sklyanin's
calculation of the Poisson brackets between matrix entries of a
monodromy matrix $M\in G$ and the observation that this Poisson structure
has the Poisson Lie group property (\cite{Drin86} Remark 5), and such results are  important in the inverse scattering method \cite{FadTak}.
The results here are `dual' to this: a space of Stokes
matrices (i.e. the ``monodromy data'' of an irregular connection)
is identified, as a Poisson manifold, with the dual group $G^*$. 

Several maps with a similar flavour have been constructed by ad hoc/homotopy theoretic means by various authors.
In the following sections we will explain that maps with the desired properties in fact arise naturally.
(Later, when considering braid group actions, we will see further applications of the above relation between Stokes data and quantum groups.)

\ 

\noindent
{\bf Ginzburg--Weinstein isomorphisms.\ }
One can also set-up the theory of Poisson Lie groups for compact groups (cf.  Lu--Weinstein \cite{LuW-PLg}):
any compact Lie group $K$ has a natural Poisson Lie group structure and the corresponding dual group $K^*$ is isomorphic to $AN$ in the Iwasawa decomposition 
$$G=KAN$$
of the complexified group  $G=K_\IC$, so for example if $K=U_n$ is the unitary group  then $A$ is the group of diagonal 
$n\times n$ matrices with real positive diagonal entries, and $N=U_+$ is the group of upper triangular unipotent complex matrices.
Thus as a manifold $K^*$ is isomorphic to $k^*$ and so one may ask if they are actually diffeomorphic as Poisson manifolds (by construction they have the same linearised Poisson structures at the origin).
The existence of such Poisson diffeomorphisms was established by Ginzburg--Weinstein \cite{GinzW} by using an indirect homotopy argument similar to an earlier  argument of  Duistermaat/Heckman.

By using involutions one may embed $K^*$ in $G^*$ and then restrict the dual exponential map to the fixed point set of the involution, and then prove that this restriction is actually a global diffeomorphism, thereby giving a new direct construction of many Ginzburg--Weinstein isomorphisms: 

\begin{thm} (\cite{smapg} for $\GL_n(\IC)$, \cite{bafi} for other $G$)
If $A_0\in k=\Lie(K)$ then the dual exponential map $\nu_{A_0}$ restricts to a global diffeomorphism of real Poisson manifolds
$$k^* \to K^*.$$
\end{thm}

Thus there are many examples of Ginzburg--Weinstein isomorphisms ``occuring in nature''.
Apparently the proof in \cite{smapg} that such maps are surjective gives a new, topological, way to show that certain Riemann--Hilbert problems have a solution.
Also, as Ginzburg--Weinstein write: such maps are no doubt related to the existence of an isomorphism $U_q(k)\cong U(k)$ of algebras (they are not isomorphic as co-algebras however since one is cocommutative).

\ 

\noindent
{\bf Duistermaat maps and Kostant's nonlinear convexity theorem.\ }\label{ssn:Duis-Kos}
Let $\lp\subset \g=\gl_n(\IC)$ denote the Hermitian matrices and let $P\subset G$ denote the positive definite Hermitian matrices.
One may identify $\lp$ with $k^*$ to give $\lp$ a  linear Poisson structure and one may identify $P\cong K^*$ using the Iwasawa and Cartan decomposition of $G$, so that $P$ also inherits a Poisson structure (cf. \cite{LuRatiu}).
The moment map for the action of the diagonal torus $T_K\subset K$ on $\lp$
is just the map
$$\delta :\lp\to \IR^n$$
taking the diagonal part of a Hermitian matrix.
Horn proved classically that if $a\in \lp$ is a diagonal matrix and $\cO$ is its conjugacy class (under the action of $K$) then the image
$$\de(\cO)\subset \IR^n$$
is a convex polytope: it is the convex hull of the $\Sym_n$ orbit of (the eigenvalues of)  $a$.
There is a nonlinear analogue of this result which goes as follows:
consider the ``Iwasawa projection'' map 
$$\wh \delta :P\to \IR^n$$
taking $g\in P\subset G$ to $\log(a)$ where $a\in A$ is the $A$ component of the Iwasawa decompositon of $g=kan\in G=KAN$.
Any conjugation orbit $\cC$ in $P$ (under $K$) is of the form
$$\cC = \exp(\cO)$$
for some orbit $\cO\subset \lp$. 
Kostant's nonlinear convexity theorem \cite{Kost73} says that 
$$\wh \de (\cC) = \de(\cO)$$
i.e. that the image of $\cC$ under the Iwasawa projection is not only a convex polytope, but that it is the same polytope as arose from the linear convexity theorem.

Now the question that Duistermaat \cite{Dui84} studied was the existence of a map $\eta$ making the following diagram commute, and in particular explaining why Kostant's nonlinear convexity result holds. 

\begin{equation}	\label{cd: texp}
\begin{array}{ccc}
  \lp & \mapright{\delta} & \IR^n \\
\mapdown{\eta} && \!\bigl\vert\!\bigl\vert \\
  P & \mapright{\wh \delta} & \IR^n.  
\end{array}
\end{equation}
Clearly taking $\eta(X)=e^X$ maps the orbits correctly, but then the diagram
does not commute. 
However one may `twist' the exponential map appropriately:

\begin{thm}[Duistermaat \cite{Dui84}] \label{thm: duistermaat} \ 

\noindent
There is a real analytic map $\psi:\lp\to K$ such that 
if one takes 
$$\eta(X)=\psi(X)^{-1}\cdot\exp(X)\cdot \psi(X)$$
then the above diagram commutes, and moreover this really is a reparameterisation, i.e.
The map $\phi_X: k\mapsto k\cdot \psi(k^{-1}Xk)$ is a diffeomorphism from
   $K$ onto $K$.
\end{thm}

Such a map clearly reduces Kostant's nonlinear convexity theorem to the linear case.
Duistermaat's motivation (also mentioned by Kostant) was to reparameterise certain integrals,
converting terms involving $\wh \delta$ into terms involving 
the linear map $\delta$.
The proof of the existence of such maps $\psi$ in \cite{Dui84} 
 involves an indirect homotopy argument.
By considering the full monodromy and Stokes data of the global connections \eqref{eq: gconn} yields a new proof of Duistermaat's theorem (in \cite{smapg} Theorem 6), and shows where such maps occur in nature, in the irregular Riemann--Hilbert correspondence: roughly speaking the linear convexity theorem considers the residue at $\infty$ and the nonlinear convexity theorem considers the Stokes data at zero, and the map $\psi$ arises as the monodromy/connection matrix relating horizontal solutions at $\infty$  to those at $0$.

\begin{figure}[h] 	
\centerline{\input{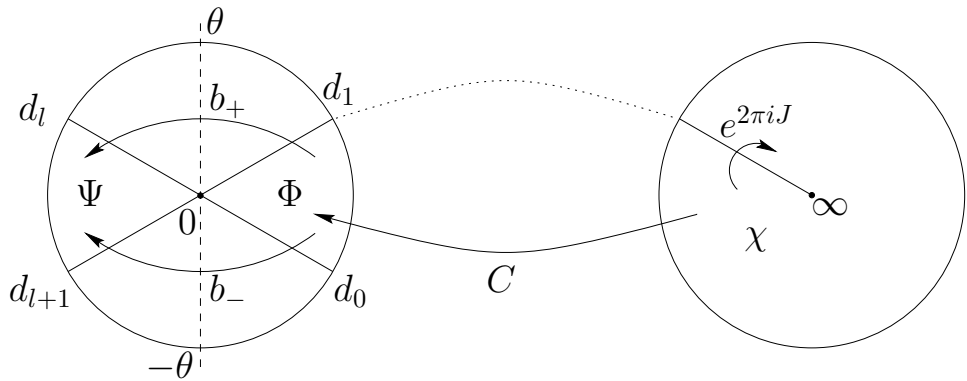}}	
\caption{Configuration in $\IP^1$ from \cite{smapg}; $\psi=C^{-1}$ is a Duistermaat map.} \label{smapg fig}
\end{figure}

In fact Kostant and Duistermaat worked with arbitrary semisimple groups (with finite centre) and our approach extends immediately to the case of complex semisimple groups (once the notion of $G$-valued Stokes data is defined, as in \cite{bafi}, and below).

\newpage

\noindent{\bf Example $\bold G$-valued Stokes data.}
We will describe the combinatorics of $G$-valued Stokes data (for $G$ a complex reductive group) in a simple example, relevant to the dual exponential map.
This example is from \cite{bafi} and extends some $\GL_n(\IC)$ cases of 
\cite{BJL79, L-R94}
(a more general situation is described in a slightly different way in \cite{gbs}).
This is the key step to defining the dual exponential map for such groups, and thus extending all the above results beyond the $\GL_n(\IC)$ case. 
When braid groups are considered in \S7 below, this naturally brings the $G$-braid groups into play, whereas the theory of isomonodromic deformations of \cite{JMU81} only involves products of type $A$ braid groups.

\begin{figure}[h] 	
\centerline{\input{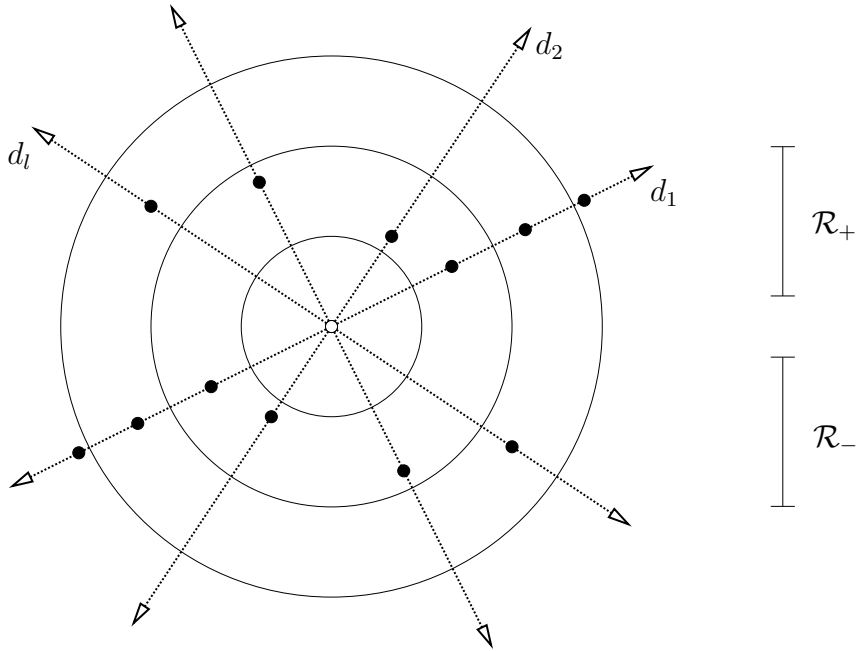}}	
\caption{Projecting the roots $\cR\subset \lt^*$ to the plane via $A_0\in \lt_\reg$, 
to define the singular directions and the Stokes groups.} \label{fig: sing dirs}
\end{figure}

Let $G$ be a connected complex reductive Lie group with maximal torus $T$ and denote the Lie algebras $\lt \subset \g$. Decomposing $\g$ with respect to $\lt$ gives the root space decomposition 
$$\g = \lt \oplus \bigoplus_{\al \in \cR} \g_\al$$
where $\cR\subset \lt^*$ is the set of roots, and for $\al\in \cR$
$$\g_\al = \{ Y\in \g \st [X,Y] = \al(X)Y \ \text{ for all $X\in \lt$}\}\subset \g$$
is the corresponding root space, which is a one-dimensional complex vector space.
Now let $\lt_\reg\subset \lt$ denote the complement of all of the root hyperplanes, the set of regular elements of $\lt$.
Choose an element $A_0\in \lt_\reg$, so that $\al(A_0)\neq 0$ for all roots $\al$.
Consider connections on the trivial principal $G$-bundle over the unit disk of the form
\beq\label{eq: conns A}
A = \left(\frac{A_0}{z^2} + \frac{B}{z} + holomorphic\right)dz
\eeq
as considered earlier.
In effect we have fixed the irregular type $Q = -A_0/z$ and are considering connections of the form $dQ+\text{\it less  singular  terms}$.
The irregular type $Q$ determines the following data:

1) a finite set $\IA\subset S^1$ of {\em singular directions} (or anti-Stokes directions) emanating from the 
singular point $z=0$ in the complex disk. These are the real directions from $0$ to the points
$$\langle A_0, \cR\rangle \subset \IC^*$$
obtained by projecting the roots onto the complex plane via the element $A_0$.
Each direction $d\in \IA$ is thus {\em supported} by some roots $\cR(d)\subset \cR$, i.e. $\cR(d)$ is the set of roots which are projected onto $d$.
See Figure \ref{fig: sing dirs}.

2) For each $d\in \IA$, a unipotent subgroup $\ISto_d\subset G$ normalized by $T$, defined as 
$$\ISto_d = \Prod_{\al\in \cR(d)} U_\al\subset G$$
where $U_\al = \exp(\g_\al)\subset G$ is the (one-dimensional) root group 
determined by $\al$, and the product may be taken in any order.
We call these the Stokes groups, and define the {\em space of Stokes data} to be 
$$\ISto(Q) = \Prod_{d\in \IA}\ISto_d$$
(here we do not take the product in $G$).
As a variety $\ISto(Q)$ is algebraically isomorphic to an affine space of dimension $\#\cR = \dim(G)-\dim(T)$.

The {\em irregular Riemann--Hilbert correspondence} associates a point of $\ISto(Q)$ to any such connection $A$. 
This lead to a natural bijection
$$\{\text{Connections $A$ in \eqref{eq: conns A}}\}/\cG_1 \ \cong \ \lt\times \ISto(Q)$$
where $\cG_1$ is the group of holomorphic maps from the unit disk to $G$ taking the value $1$ at $z=0$. This statement is equivalent to (the $k=2$ case of) \cite{bafi} Theorem 2.8. More recently  (cf. \cite{gbs} v3, Appendix A) such statements may be ``upgraded''
to an equivalence of categories between  connections with fixed irregular types and {\em Stokes $G$-local systems}.

Now we will sketch how to associate Stokes data to a connection $A$.
In effect we prove the 
multisummation approach (of \cite{ramis-factn, BBRS91, MR91, L-R94} etc.) goes through.
(In the present example Borel summation is sufficient though.) 
The crucial fact (see \cite{bafi} Lemma A.5) is that complex reductive Lie groups 
(those whose representations are completely reducible) are affine algebraic groups, and so the fact that multisummation is a morphism of differential algebras implies things work nicely. (This also implies there is no trouble extending this approach to any affine algebraic group.)
The basic statements are as follows:

1) (\cite{bafi} Lemma 2.1.) There is a unique element $\La\in \lt$ and formal gauge transformation $\wh F\in G\flb z \frb$ with $\wh F(0)=1$ such that 
$$A= \wh F [ A^0], \qquad\text{where}\qquad A^0 :=\left(\frac{A_0}{z^2} + \frac{\La}{z}\right)dz.$$

\noindent
Thus $A$ is {\em formally} isomorphic to an abelian ($\lt$-valued) connection $A^0$.

2) (\cite{bafi} Theorem 2.5.) If $\Sect_i$ is a sector (in the unit disc) bounded by two consecutive singular directions, then there is a preferred choice 
$\Si_i(\wh F):\Sect_i\to G$ 
of an analytic isomorphism between $A$ and $A^0$ asymptotic to $\wh F$ at $z=0$.

3) If $\Sect_i, \Sect_{i+1}$ are consecutive sectors, abutting at $d\in \IA$, then
the fundamental solutions 
$$\Phi_i = \Si_i(\wh F) z^\La e^Q:\Sect_i\to G,\qquad 
  \Phi_{i+1} = \Si_{i+1}(\wh F) z^\La e^Q:\Sect_{i+1}\to G$$
of $A$ may be analytically continued across $d$ and then 
$$\Phi_i = \Phi_{i+1}\circ K_d$$ 
for a unique ($z$-independent) element $K_d\in \ISto_d$ (\cite{bafi} Lemma 2.7).
(Here $z^\La e^Q$ denotes a fundamental solution of $A^0$, continuous across $d$.)

4) Repeating for each $d\in \IA$ 
yields a {\em surjective} map 
$$\{\text{Connections $A$ in \eqref{eq: conns A}}\} \ \to \ \lt\times \ISto(Q)$$
taking the Stokes data $K_d\in \ISto_d$ and the ``exponent of formal monodromy'' 
$\La\in \lt$. The fibres of this map are precisely the $\cG_1$ orbits (cf. \cite{bafi} Theorem 2.8).

5) Finally we can reorganise the Stokes data. If we choose a sector $\Sect_0$ (bounded by consecutive singular directions) and let $\Sect_l = -\Sect_0$ be the opposite sector, then the singular directions 
$d_1,\ldots,d_l$ one crosses on going from $\Sect_0$ to $\Sect_l$ in a positive sense, support a system of positive roots $\cR_+=\cR(d_1)\cup\cdots\cup \cR(d_l)\subset \cR$.
The product (in $G$) of the corresponding Stokes groups is isomorphic (as a space) to the unipotent radical $U_+$ of the Borel subgroup $B_+$ determined by $\cR_+$
(\cite{bafi} Lemma 2.4). Similarly going from $\Sect_l$ to $\Sect_0$ in a positive sense
yields the unipotent radical $U_-$ of the opposite Borel.
In this way the choice of $\Sect_0$ determines an isomorphism $\ISto(Q)\cong U_+\times U_-$, and in turn, adding in $\La$, we obtain
$$\lt\times\ISto(Q)\cong \lt\times U_+\times U_-\cong G^*$$
which, as a space, is the simply connected dual Poisson Lie group $G^*$.

\newpage\noindent%
{\large \bf5. Fission}\label{sn: fission}

The articles \cite{fission, gbs} introduce a new operation ``fission'', complementary (and not inverse to) the ``fusion'' operation of Alekseev et al \cite{AMM} that they used to construct symplectic moduli spaces of flat connections on Riemann surfaces.
In brief, when one considers symplectic moduli spaces of meromorphic connections on Riemann surfaces, fusion enables 
induction with respect to the genus and number of poles, whereas fission 
enables induction with respect to the {\em order} of the poles.

\noindent{\bf Background.\ }
The quasi-Hamiltonian approach \cite{AMM} to building symplectic moduli spaces of 
representations  of the fundamental group of
a Riemann surface involves starting with some simple pieces, and then using two operations: {\em fusion} and {\em reduction}.
(See e.g. \cite{Mein12} for a recent introduction to these ideas.)
Both of these operations are (quasi-)classical analogues of operations in conformal field theory, related to gluing together Riemann surfaces with boundary.
(A half-way step between the physics and the algebraic approach of \cite{AMM} are  the 
Hamiltonian loop-group spaces of Donaldson \cite{Don-bvym}  and  Meinrenken--Woodward \cite{MeiWoo-fusion}---fusion is described at this level in \cite{MeiWoo-fusion} \S4.1).

Fusion puts a ring structure on the category of quasi-Hamiltonian $G$-spaces:
One may attach a quasi-Hamiltonian $G$-space $M(\Si)=\Hom(\pi_1(\Si),G)$ to any Riemann surface $\Si$ with exactly one boundary component and the fusion product of two such spaces $M(\Si_1), M(\Si_2)$ 
is the space attached to the surface $\Si_3$ obtained by gluing 
$\Si_1$ and $\Si_2$ into two of the holes of a three-holed sphere:
 $$M(\Si_1) \fus M(\Si_2) \ = \ M(\Si_3)$$

\begin{figure}[h]
\begin{center}
\includegraphics[angle=270, width=11cm, trim=0 20 0 -50, clip]{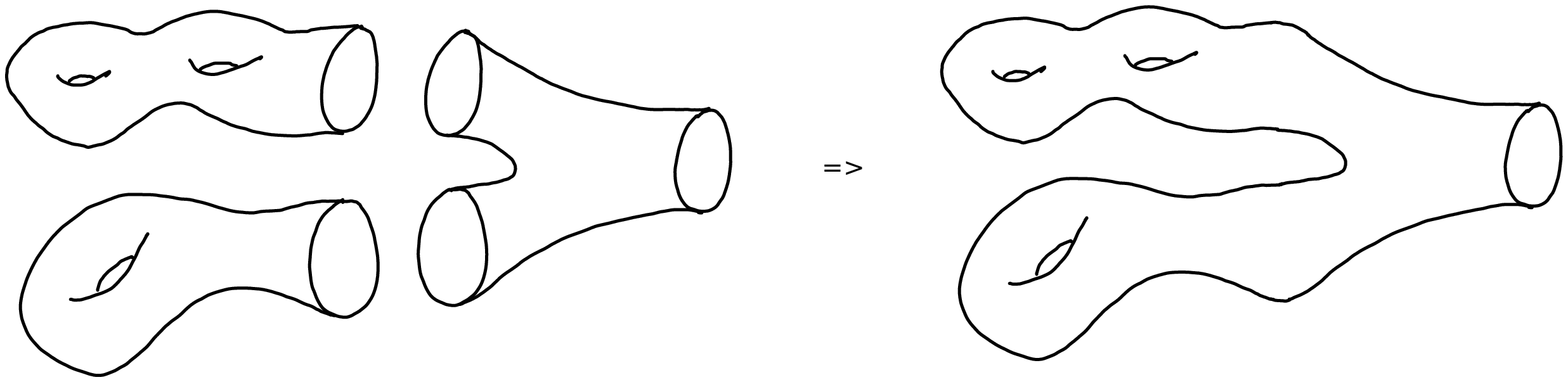}
\caption{Fusion}
\end{center}
\end{figure}

\noindent{\bf Main results.}
A new sequence of quasi-Hamiltonian spaces was constructed in 
\cite{saqh, fission, gbs} and these may be used to replace the three-holed sphere
in the definition of the fusion operation above, yielding a sequence of new operations, which we call ``fission''.
They enable one to combine quasi-Hamiltonian spaces for different structure groups $G$.

The new spaces are as follows. Let $G$ be a connected complex reductive group, $P_+,P_-\subset G$ a pair of opposite parabolic subgroups,
let $H=P_+\cap P_-$ be their common Levi factor and let $U_\pm\subset P_\pm$ be their unipotent radicals.

\begin{thm}[\cite{saqh, fission, gbs}]
For any integer $r\ge 1$ the space
$$\gahr := G\times (U_+\times U_-)^r \times H$$
is a quasi-Hamiltonian $G\times H$-space.
\end{thm}

More precisely this is proved in \cite{saqh} in the case when $P_\pm$ are opposite Borels (so that $H$ is a maximal torus), in \cite{fission} for $r=1$ (with $H$ possibly non-abelian) and in general in \cite{gbs}.
This gives a way to break the structure group $G$ to the subgroup $H$ (the Dynkin diagram of $H$ is obtained from that of $G$ by deleting some nodes)---this motivates the name ``fission''.
Some quasi-Hamiltonian spaces of M. Van den Bergh \cite{vdb-doublepoisson, vdb-ncqh, yamakawa-mpa} arise from simple examples of these fission spaces (\cite{gbs} \S4).

One may see these spaces yield new operations as follows.
First, the usual fusion picture above may be rephrased as follows. 
Let $S$ denote the three-holed sphere. 
This yields a quasi-Hamiltonian $G^3$-space 
$M(S) = \Hom(\Pi_1(S),G)$, where $\Pi_1(S)$ denotes the fundamental groupoid of $S$ with one basepoint on each boundary component (cf. \cite{gbs} Theorem 2.5).
Then fusion amounts to the following {\em gluing}:
$$M(\Si_1)\fus M(\Si_2) = M(\Si_1)\glue{G} M(S) \glue{G} M(\Si_2)$$

\noindent
where the symbol $\glue{}$ denotes the gluing (cf. \cite{fission} \S5).  
Since $M(S)$ is a quasi-Hamiltonian $G^3$-space, and each gluing absorbs a factor of $G$, the result is a quasi-Hamiltonian $G$-space, as expected.

Now, for the fission spaces, typically $H$ will factor as a product of groups (e.g. if $G=\GL_n(\IC)$ then $H$ is a ``block diagonal'' subgroup). 
Suppose $H=H_1\times H_2$ for definiteness (the generalisation to arbitrarily many factors is immediate).
Thus $\gahr$ is a quasi-Hamiltonian $G\times H_1\times H_2$ space.
For example this enables us to construct a quasi-Hamiltonian $G$-space 
$$M_1\glue{H_1} \gahr \glue{H_2} M_2$$

\noindent
for any integer $r$, out of quasi-Hamiltonian $H_i$-spaces $M_i$ ($i=1,2$), e.g. we could take  $M_i = \Hom(\pi_1(\Si_i),H_i)$.
Thus the fission spaces yield many new operations on the category of quasi-Hamiltonian spaces (without fixing the group $G$ beforehand).
One may picture these operations as indicated in Figure \ref{fig: Y}.

\begin{figure}[ht]
	\centering
	\input{annulus.and.Y.pstex_t}
	\caption{Fission%
}\label{fig: Y}
\end{figure}

One application of the fission operations is to construct the wild character varieties
(see p.\pageref{sn: wcvs}). Surprisingly it turns out that  many other algebraic symplectic manifolds may be constructed in this way as well, such as all of the so-called multiplicative quiver varieties (\cite{gbs} Corollary 4.3).

\newpage
\noindent%
{\bf\large6. Wild character varieties} \label{sn: wcvs}

In the articles \cite{saqh, fission, gbs} the author has constructed the wild character varieties.
These are symplectic algebraic varieties which generalise the 
complex character varieties of Riemann surfaces.
The well-known braid and mapping class group actions on the character varieties are also generalised in \cite{bafi, gbs}.

\noindent{\bf Background.}
Given a Riemann surface $\Si$ (maybe open or with boundary), 
and a Lie group $G$ many people have studied the 
{\em character variety} 
$$\MB(\Si,G)=\Hom(\pi_1(\Si),G)/G$$ 
of representations of the fundamental group of $\Si$ into $G$.
If $G$ is a complex reductive group then there is a natural holomorphic Poisson structure on $M$, first considered analytically by Atiyah--Bott
\cite{AB83}, then understood in terms of group cohomology by Goldman \cite{Gol84} and 
subsequently studied purely algebraically by many people such as 
\cite{Kar92,FR,Weinstein-floermem,Jeffrey-95,AMR,GHJW,AMM}.
See e.g. \cite{Aud95long} or \cite{sikora-charvars} for an overview.

If $\Si$ is in fact a smooth complex algebraic curve (possibly punctured)
and $G=\GL_n(\IC)$ 
then Deligne's Riemann--Hilbert correspondence \cite{Del70} implies that 
$\MB$ is isomorphic 
to the space of algebraic connections on rank $n$ vector bundles on $\Si$ with {\em regular singularities} at the punctures.
On a curve, a connection has regular singularities if and only if it may be obtained by restricting a meromorphic connection on the compact curve which only has {\em simple poles} at the marked points.

This raises the question of constructing the analogue of the spaces $\MB$ which classify more general connections on curves, not satisfying this regularity assumption, the {\em irregular} connections.
The irregular Riemann--Hilbert correspondence was first considered by Birkhoff \cite{birkhoff-1913} but was only fully worked out on curves 
quite recently (see \cite{malg-book} and references therein, and the recently published letters of Deligne--Malgrange--Ramis \cite{DMR-ci}).
In brief one adds some ``Stokes data'' at each marked point, and there are various ways to package this extra data. 
None of this work considers moduli spaces or symplectic structures however.

One motivation for pursuing this is that the wild character varieties admit interesting braid group actions, the simplest case of which is known
(\cite{bafi}) to underly the so-called quantum Weyl group actions.

The only previous study of {\em moduli spaces} of Stokes and monodromy  data in any serious generality 
is in the integrable systems literature: 
Jimbo--Miwa--Ueno \cite{JMU81} considered 
the case of certain connections on the Riemann sphere with 
just one level.
On the other hand Flaschka--Newell \cite{FN82} considered symplectic structures in some $\GL_2(\IC)$ examples.

\noindent{\bf Main results.}
The main results are stated succinctly in the introduction to \cite{gbs}.
In brief the (complexification of the) quasi-Hamiltonian approach of 
Alekseev--Malkin--Meinrenken \cite{AMM} to $\MB$ says that $\MB$ arises as a finite dimensional algebraic {\em multiplicative}  
symplectic quotient of a smooth (finite dimensional) affine variety. 
This gives a purely algebraic approach to the Poisson structure on $\MB$.

The articles  \cite{saqh, fission, gbs} show how to extend this approach to the irregular case (the articles are in increasing generality, leading up to the case of any connected complex reductive group $G$, with meromorphic connections having arbitrary unramified formal normal forms on arbitrary genus smooth algebraic curves).

To describe this
it is convenient to define the notion of an ``irregular curve'' $\Si$ 
to be a smooth complex algebraic curve together with some marked points $a_1,\ldots,a_m\in \Si$ plus the extra data of an {\em irregular type} 
$Q_i$  at each marked point.
It $z$ is a local coordinate vanishing at $a_i$ then
$Q_i = A_r/z^r+\cdots A_1/z$ for elements $A_i\in \lt$ in a Cartan subalgebra of $\g=\Lie(G)$.
This generalises the notion of a curve with marked points, and it turns out  to be very useful to view the irregular types as analogous to the moduli of the curve in this way: they behave just like the moduli of the curve (and similarly lead to interesting braid group actions when varied).
Given an irregular type $Q_i$ at $a_i$ 
we consider connections locally of the form 
$$dQ_i + \text{ less singular terms}$$
near $a_i$, so that $Q_i$ determines the irregular part of the connection, and such connections have solutions involving essentially singular terms of the form $e^{Q_i}$. 
Then for any irregular curve $\Si$, \cite{gbs} defines a certain groupoid $\Pi$ and the space 
$\Hom_\IS(\Pi,G)$ of Stokes representation of $\Pi$.
This has a natural action of the group 
$\bH:=H_1\times\cdots \times H_m$ where $H_i = C_G(Q_i)$ is the centraliser of $Q_i$.
The main result is then:

\begin{thm}[\cite{gbs}]
The space $\Hom_\IS(\Pi,G)$ of Stokes representations is a smooth affine variety and is a quasi-Hamiltonian $\bH$-space, 
where $\bH=H_1\times \cdots \times H_m\subset G^m$.
\end{thm}

This implies that the quotient $\Hom_\IS(\Pi,G)/\bH$ (the wild character variety), which classifies
meromorphic connections with the given irregular types, inherits an algebraic Poisson structure. 
Its symplectic leaves are obtained by fixing a conjugacy class 
$\cC_i\subset H_i$ for each $i=1,\ldots,m$.
In the regular singular case, when each irregular type $Q_i=0$, the space
$\Hom_\IS(\Pi,G)$ is just the space of all representations of the fundamental groupoid of $\Si\setminus\{a_i\}$ (with a basepoint near each puncture) in the group $G$, and 
$\bH=G^m$ so that $\Hom_\IS(\Pi,G)/\bH \cong \MB(\Si\setminus \{a_i\},G)$ and we recover the original picture.

The article \cite{gbs} also characterises the stable points of $\Hom_\IS(\Pi,G)$ in the sense of geometric invariant theory (for the action of $\bH$), using the Hilbert--Mumford criterion, shows there are lots of examples when the quotients are well-behaved
and describes the irregular analogue of the Deligne--Simpson problem.

As a corollary of this approach,
in \cite{gbs} Corollary 9.9 it is proved that, with fixed generic conjugacy classes $\cC_i\subset H_i$,
the wild character varieties are {\em smooth} symplectic affine varieties, 
in the case $G=\GL_n(\IC)$.
(This gives a direct algebraic description of the spaces underlying the \hk manifolds of \cite{wnabh}.)
This result alone probably justifies the quasi-Hamiltonian approach---it generalises a result of Gunning \cite{gunning-lovb} \S9 in the holomorphic case $m=0$, obtained by explicitly differentiating the monodromy relation, something that seems daunting in the present set-up.

\newpage

\noindent%
{\bf \large7. Braid group actions from isomonodromy}\label{sn: bafi}

In the articles \cite{bafi, gbs} the author defined the notion of $G$-valued Stokes data (for $G$ a connected complex reductive group) and set up the theory of isomonodromic deformations in this context
generalising some work of Jimbo--Miwa--Ueno \cite{JMU81}.
Geometrically this amounts to defining the notion of an admissible family of irregular curves and showing that a natural nonlinear flat connection exists on 
the bundle of wild character varieties associated to an admissible family of irregular curves \cite{gbs}.
The monodromy of these nonlinear connections was computed explicitly in some cases in \cite{bafi} and shown to yield the $G$-braid group actions underlying those of the so-called quantum Weyl Group.

\noindent
{\bf Background.} \ 
The classical theory of monodromy preserving deformations of linear differential equations on the Riemann sphere (or ``isomonodromic deformations'') 
was revisited and extended in the 
early 1980's by Jimbo--Miwa--Ueno \cite{JMU81, JM81} (see also Flaschka--Newell \cite{FN80}) due to the appearance of such deformations in physical problems (cf. \cite{JMMS, WMTB, BMW, myers}).
Later in \cite{smid} the author revisited \cite{JMU81, JM81} from a (symplectic) geometric perspective and rephrased some of their results in terms of nonlinear connections on fibre bundles.
(As mentioned in \cite{smid} this was motivated by the appearance of certain examples of such isomonodromic deformations in the classification of two dimensional topological quantum field theories/Frobenius manifolds.)
The parallel with Simpson's Gauss--Manin connection in nonabelian cohomology
\cite{Sim94ab} was also noted (\cite{smid} introduction and  \S7). 
In this work Simpson shows there is a natural flat nonlinear connection 
on the bundle of first nonabelian cohomologies associated to any family of smooth projective varieties. 
But the De\! Rham description of the first nonabelian cohomology is as the moduli space of flat holomorphic connections on vector bundles, and so we see the nonlinear connections of \cite{JMU81, JM81} (as described in \cite{smid}) are analogues of this when one extends from holomorphic to meromorphic connections (and takes the underlying projective variety to be the Riemann sphere).

A crucial difference however is that Jimbo--Miwa--Ueno understood that in the case of irregular meromorphic connections there are many more independent deformation parameters beyond the moduli of the underlying Riemann sphere with marked points. These extra parameters (the ``irregular times'') control the irregular type of the connections.
Thus the picture of  \cite{JMU81, JM81} is not just the extension of the nonabelian Gauss--Manin connection to the case of certain quasi-projective varieties, but an extension involving new deformation parameters, hidden in the classical algebro-geometric picture of deriving flat connections from families of varieties.

\noindent{\bf Main Results.}
Since flat connections are not easy to come by in mathematics (and intrinsic geometric ones especially so) the articles \cite{bafi, gbs} %
further pursued and generalised the irregular times of Jimbo et al, and the resulting 
theory of isomonodromic deformations (irregular nonabelian Gauss--Manin connections).
The author's feeling is that these extra parameters should be taken as seriously as the moduli of the underlying Riemann surface with marked points.

To this end, given a complex reductive group $G$  with a fixed maximal torus, the article  \cite{gbs} defines the notion of an ``irregular curve'' $\Si$ consisting of a compact smooth complex algebraic curve, plus some marked points and an irregular type at each marked point.
If the irregular types are zero this specialises to the notion of curve with marked points.
As explained in \S6 above, to any irregular curve there is a canonically associated wild character variety $\MB(\Si) = \Hom_\IS(\Pi,G)/\bH,$ which is naturally a Poisson variety.
Then \S10 of \cite{gbs} defines the notion of an ``admissible family'' of irregular curves, generalising the notion of deforming a smooth curve with marked points such that it remains smooth and none of the points coalesce.  

Given an admissible family $\pi:\Si\to \IB$ of irregular curves over a base space $\IB$ then one can consider the family of wild character varieties  
$\MB(\Si_p)$ as $p\in \IB$ varies, where $\Si_p$ is the irregular curve 
$\pi^{-1}(p)$ over $p\in \IB$.

\begin{thm}(\cite{gbs} \S10) \label{thm: loc syst of vars}
The varieties $\MB(\Si_p)$ assemble into a local system of Poisson varieties over $\IB$.
\end{thm}
This means that there is a fibre bundle $\pr:M\to \IB$ 
such that $\pr^{-1}(p) = \MB(\Si_p)$ for any $p\in \IB$, with a complete flat Ehresmann connection on it (the irregular isomonodromy connection); for any points $p,q\in \IB$ and path $\ga$ in $\IB$ from $p$ to $q$, there is a canonical algebraic Poisson isomorphism $\MB(\Si_p)\cong \MB(\Si_q)$, only dependent on the homotopy class of $\ga$.
Consequently there is an algebraic Poisson action of the fundamental group $\pi_1(\IB,p)$ on $\MB(\Si_p)$ for any basepoint $p\in \IB$. 
(This generalises the well known braid and mapping class group actions in the ``usual'' theory.)

This extends the viewpoint of Jimbo et al \cite{JMU81, JM81} in several ways, since we allow: 
1) $G$ to be any complex reductive group, 2) $\Si$ to have any genus, 
3) any unramified irregular types (e.g. the leading coefficients may have repeated eigenvalues in the general linear case), and also 4) since we consider algebraic Poisson/symplectic structures and show they are preserved.
In fact the innovation of allowing  $G$ to be any complex reductive group, and phrasing the Stokes data in terms of the roots, appeared in the earlier article \cite{bafi}. An application of this will be described in the following subsection.

Note that from the viewpoint we started with in \cite{Boa, smid} (due to \cite{cec-vafa-nequals2classn, Dub95long-with.ApE.note}) of Stokes data classifying topological quantum field theories, with matrix entries counting BPS states (or solitons) going between $n$ vacua, the idea (of \cite{bafi}) of passing from $\GL_n(\IC)$ to another algebraic group is quite bizarre since it would correspond to passing to a $\g=\Lie(G)$-valued Cartan matrix (cf. e.g. \cite{cec-vafa-nequals2classn} (6.21), \S7.1, \cite{Dub95long-with.ApE.note} (H26), pp.263-4).
Nonetheless this idea is used in some recent work on wall crossing of BPS states (see e.g. \cite{btl1}).

\noindent{\bf Geometric origins of the quantum Weyl group.\ }
Recall from \S4 we have a new 
geometric/moduli-theoretic viewpoint on the theory of quantum groups: the Drinfeld--Jimbo quantum group is the quantisation of a very simple moduli space of irregular connections having a pole of order two.
Thus one would expect other features of the quantum group to appear geometrically as well.
The so-called ``quantum Weyl group'' action is essentially an action of the $G$-braid group on the quantum group and was defined explicitly by  
 Lusztig \cite{Lus90b}, Soibelman \cite{Soib} and Kirillov--Reshetikhin \cite{KResh}, via generators and relations.
The quasi-classical limit of this action, a Poisson action of the $G$-braid group on $G^*$, was computed explicitly by 
De\! Concini--Kac--Procesi \cite{DKP}.
On the other hand since we have identified $G^*$ with a space of Stokes data (i.e. it is essentially a wild Betti space $\cM_B$), as above, by integrating the isomonodromy connection, we obtain a nonlinear discrete group action on $G^*$, analogous to the usual mapping class group action on the character varieties. 
In the present context the space $\IB$ of deformations is the 
space of regular elements $A_0\in \lt_\reg$, whose fundamental group is the pure $G$-braid group. This may be extended to the full braid group (adding in the finite Weyl group) to obtain the following statement:

\begin{thm}(\cite{bafi} Theorem 3.6.)
The De\! Concini--Kac--Procesi action of the $G$-braid group on $G^*$ coincides with the isomonodromy action, and so the quantum Weyl group action quantises the isomonodromy action.
\end{thm}

Some further aspects of this story are also elucidated in \cite{bafi}, going around the square in Figure \ref{fig: hopf algebras} on p.\pageref{fig: hopf algebras} above.
The quantum Weyl group action is defined at the top, on $\U_q\g$, and De\! Concini--Kac--Procesi followed the arrow 1 down to the left, and the above theorem shows the braid group action on $G^*$ they computed in this way comes from isomonodromy. 
However the isomonodromy action is obtained by integrating a nonlinear connection, and this connection is essentially equivalent to the explicit system of nonlinear differential equations 
\beq\label{eq: nldif}
dB = \left[B, \ad^{-1}_{A_0} [dA_0,B]\right]
\eeq
for $B\in \g\cong \g^*$ as a function of $A_0\in\lt_\reg=\IB$ 
(cf. \cite{bafi} (4.3)).
This arises by passing to the other side of the irregular Riemann--Hilbert correspondence, essentially conjugating by  the dual exponential maps 
$\nu_{A_0}:\g^*\to G^*$ (i.e. following arrow 3 down to $\g^*$).
In other words this differential equation is the infinitesimal manifestation of the braid group action at the level of $\g^*$.
Finally we can ascend the arrow 4: This arrow corresponds to the PBW quantisation of $\Sym\g$ into $U\g$; we apply the symmetrisation map 
to the Hamiltonians for the system \eqref{eq: nldif} to obtain a
a flat connection on the trivial $U\g$ bundle over $\IB=\lt_\reg$ (as written in \cite{bafi} Proposition 4.4 and (4.7)).
This connection is the simplest irregular analogue of the Knizhnik--Zamolodchikov (KZ) connection and was guessed, and shown to be flat, by
De\! Concini (unpublished), and Millson--Toledano Laredo \cite{millson-toledano, VTL-duke}\footnote{This connection was called the DMT connection in \cite{bafi}, but in fact a similar, indeed slightly more general, connection appeared before \cite{millson-toledano, VTL-duke} in work of Felder et al  \cite{FMTV}
(this reference should have been included in \cite{bafi}, and apologies are due to the authors of \cite{FMTV}).} (who had no idea that their work was related to irregular connections).   We thus have a simple intrinsic {\em derivation} of this irregular KZ connection from the isomonodromy Hamiltonians (and have connected all the vertices of the square in Figure \ref{fig: hopf algebras}).

\ppb{
\newpage

\noindent$\bullet$
{\bf Irregular mapping class groups.\ }
The article \cite{gbs} also introduces the notion of an ``irregular curve'' as a Riemann surface with some marked points and some extra data (irregular types) at the marked points, and shows that a wild Betti space may be canonically associated to any irregular curve, and that they  behave well under `admissible' deformations of the irregular curve: the wild Betti spaces form a ``Poisson local system of varieties'', \cite{gbs} \S10.
In other words, if we have a space $\IB$ parameterising an admissible family of irregular curves,   the corresponding family of wild Betti spaces fit together into the fibres of a fibre bundle over $\IB$  with a flat, algebraically integrable connection, preserving the natural Poisson structures.
Such deformations generalise the notion of deforming a curve with marked points such that it remains smooth and the marked points do not coalesce: The fundamental group of the space of admissible deformations is (thus) a natural generalisation of the mapping class groups/braid groups of Riemann surfaces. And the fact that the wild Betti spaces  
form a Poisson local system of varieties, means that such fundamental groups (the wild/irregular mapping class groups) 
act naturally on the wild Betti spaces preserving their symplectic/Poisson structures (generalising the well-known mapping class group actions on the character varieties).
The simplest example of such wild mapping class group actions is the quasi-classical limit of the quantum Weyl group, as shown in \cite{bafi}.

}

\newpage
\noindent%
{\bf\large8. Logahoric connections on parahoric bundles} \label{sn: logahoric}

The article \cite{logahoric} defines the notion of logahoric connections (i.e. an analogue of logarithmic connections on parahoric bundles) and shows that there is a  Riemann--Hilbert correspondence for them.
The corresponding (local) monodromy/Betti data consists of pairs $(M,P)$ where 
$M\in G$ is the local monodromy, $P\subset G$ is a (weighted) parabolic subgroup and $M\in P$.   
If $P$ is a Borel subgroup then such data appears in the multiplicative Brieskorn--Grothedieck--Springer resolution; 
we thus obtain a moduli-theoretic realisation of the multiplicative Brieskorn--Grothedieck--Springer resolution.
We also construct the natural multiplicative symplectic structures (i.e. quasi-Hamiltonian structures) so the resolution map is now the {\em group valued} moment map.

\noindent{\bf Background.}
In his work on the nonabelian Hodge correspondence on noncompact curves, 
Simpson \cite{Sim-hboncc} established a Riemann--Hilbert correspondence for 
``tame filtered $\cD$-modules'' on a curve. 
These objects may be understood as {\em logarithmic connections on parabolic vector bundles}.
Recall a parabolic vector bundle \cite{MehSes} on a smooth compact curve $\Si$  with marked points $a_1,\ldots, a_m\in \Si$ consists
of a holomorphic vector bundle $V\to \Si$ together with a filtration in the fibre $V_{a_i}$ of $V$ at each marked point.
A filtration consists of a weighted flag (traditionally the weights are rational numbers in $[0,1)$, but for the full nonabelian Hodge correspondence it is necessary to work with all the real numbers in this interval). 
On a curve a {\em logarithmic connection} is just a meromorphic connection having poles of order $\le 1$. On a parabolic vector bundle, the residues of the connections should preserve the flags.

Simpson sets up a correspondence between these objects and {\em filtered local systems} on the punctured curve: this amounts to a representation of the fundamental group plus, near each puncture, a filtration in a nearby fibre preserved by the local monodromy (now the weights are arbitrary real numbers, not restricted to be in an interval).
The filtration encodes the growth rate of solutions 
(and is closely related to 
the $\IZ$-filtrations of Levelt \cite{levelt61} (2.2) in the case of logarithmic connections on usual vector bundles).

Now suppose we replace the structure group $G=\GL_n(\IC)$ used above by an arbitrary connected complex reductive group $G$.
Then it is reasonably clear how to generalise the notion of filtered local system: one takes a representation of the fundamental group of the punctured curve into $G$---i.e. a $G$-local system---plus a weighted parabolic subgroup in a fibre near each marked point, preserved by the local monodromy (we call this a ``filtered $G$-local system'' cf. \cite{logahoric} Remark 2). 

\noindent
{\bf Question.} Is there a Riemann--Hilbert correspondence for filtered $G$-local systems, and if so what are the corresponding connection-like objects?

The point is that one does {\em not} get a full correspondence by considering 
logarithmic connections on parabolic $G$-bundles. %
Said differently this question is asking: What are the basic objects that should appear in the tamely ramified nonabelian Hodge correspondence for $G$-bundles on noncompact curves?

\newpage
\noindent
{\bf Main results.}
The answer %
is to consider an analogue of logarithmic connections when one replaces
a parabolic $G$-bundle by a {\em parahoric bundle}, i.e. a torsor under a (weighted) parahoric group scheme on the compact curve.

The notion of ``quasi-parahoric bundle'' has been studied recently by various authors, 
such as \cite{pap-rap, heinloth-uni, yun2}: it is a torsor under a parahoric (Bruhat--Tits) group scheme on the compact curve.
This generalises the notion of quasi-parabolic vector bundle due to Mehta--Seshadri 
\cite{MehSes} (this is a parabolic bundle when one forgets the weights), and the 
notion of quasi-parabolic $G$-bundles (considered e.g. in \cite{lasz-sorg}). 
In general one does not have an underlying $G$-bundle on the compact curve.

Thus the first step is to define the notion of ``weight'' for a quasi-parahoric bundle:
in brief near each marked point a quasi-parahoric bundle amounts to the choice of a parahoric subgroup of the local loop group, and these are classified by the facettes in the Bruhat--Tits building.
But the  Bruhat--Tits building \cite{BrTits-I} p.170 is built out of real vector spaces, the apartments (although it is often viewed as a simplicial complex).
So we can define a {\em weighted parahoric subgroup} to be a 
{\em point} of the  Bruhat--Tits building (\cite{logahoric} Defintion 1, p.46).
This yields the notion of parahoric bundle.

Next we need to find the right notion of singular connection. 
Locally a $G$-bundle corresponds to the parahoric subgroup $G\flb z \frb \subset G\flp z \frp$, and a logarithmic connection is a connection having a pole of order one, i.e. is represented by an element of $\g\flb z \frb dz/z$, where $\g =\Lie(G)$.
In general we {\em define} a ``logahoric'' connection (or a ``tame parahoric connection'') to be a connection with a pole of order one more than that permitted by the parahoric level structure, \cite{logahoric} \S3.
(In general such connections may have arbitrary order poles, but will always be regular singular connections.)
The main result is that there is a precise correspondence between these objects and filtered $G$-local systems. This follows from the local correspondence which may be stated as follows.

\begin{thm}[\cite{logahoric}]
There is a canonical bijection between $LG=G\flp z\frp$ orbits of tame parahoric connections and $G$ orbits of enriched monodromy data:
$$\Bigl\{ (A,p) \st p\in \cB(LG), A\in \cA_p\Bigr\}/LG \quad \cong\quad 
 \Bigl\{ (M,b) \st b\in \IB(G), M\in P_b\Bigr\}/G.$$
\end{thm}
Here $\cB(LG)$ is the Bruhat--Tits building (the space of weighted parahoric subgroups of $LG$), and $\IB(G)$ is the space of weighted parabolic subgroups of $G$ ($\cA_p$ is the space of logahoric connections determined by $p\in \cB(LG)$ and $P_b\subset G$ is the parabolic underlying $b\in \IB(G)$).

The relation to Grothendieck's simultaneous resolution, and what we thus learn about its geometry,  will be described in the next two pages.

\noindent{\bf Further developments.}
Six months later (after posting on the arXiv and submitting \cite{logahoric} 
to a journal with Seshadri on the editorial board) 
Balaji--Seshadri \cite{BS10} used a similar (but slightly less general) notion of weights for parahoric torsors,
and they established an analogue of the Mehta--Seshadri theorem 
(although they did not consider the full Riemann--Hilbert correspondence).

\newpage
\noindent%
{\bf Geometry of the Brieskorn--Grothendieck--Springer resolution. }
First we will quote (from Brieskorn's ICM talk \cite{briesk-icm1970}) 
a result proved by Grothendieck (and already by Springer \cite{springer-unip} as far as the unipotent fibre is concerned). 
Let $G$ be a (simply-connected) semisimple complex algebraic group, with maximal torus $T$ and Weyl group $W$.
Then $T/W$ 
parameterises the conjugacy classes of semisimple elements of $G$; there is a map $\psi:G\to T/W$ taking an element to the class of its semisimple part (using the Jordan decomposition). The fibres of $\psi$ are unions of conjugacy classes of $G$.
Let $B_0\subset G$ be a Borel subgroup containing $T$, let 
$\cB\cong G/B_0$ denote the variety of Borel subgroups of $G$, and define
$$\wt G = \{(M,B)\in G\times \cB \st M\in B\}.$$
to be the set of pairs consisting of a group element $M$ and a Borel subgroup $B$ containing $M$.
Projection onto the first factor gives a map $\pi:\wt G \to G$ and there is a natural map $\wt \psi:\wt G\to T$ such that the following diagram commutes.

\begin{thm}
The following diagram is a simultaneous resolution of the singularities of the fibres of $\psi:G\to T/W$:
\begin{equation}	\label{cd: BGS}
\begin{array}{ccc}
  \wt G & \mapright{\pi} & G \\
\mapdown{\wt\psi} && \mapdown{\psi} \\
  T & \mapright{\pr} & T/W.  
\end{array}
\end{equation}
\end{thm}

\noindent
In particular for any $t\in T$ the map 
$\pi:\wt \psi^{-1}(t)\to \psi^{-1}(\pr (t))$ 
is a resolution of singularities.
The fibres $\wt \psi^{-1}(t)\subset \wt G$ are of the form 
$G\times_{B_0} tU$ where $U\subset B_0$ is the unipotent radical. 
For more details see \cite{briesk-icm1970}, or Slodowy \cite{slodowy-book} Theorem 4.4, or Steinberg \cite{steinberg-unip} \S6.

There is a similar ``additive''  statement on the Lie algebra level 
(\cite{slodowy-book} \S4.7) with $G$ replaced by $\g=\Lie(G)$ and Borel subgroups by Borel subalgebras.
\begin{equation}	\label{cd: aBGS}
\begin{array}{ccc}
  \wt \g & \mapright{\pi} & \g \\
\mapdown{\wt\psi} && \mapdown{\psi} \\
  \lt & \mapright{\pr} & \lt/W.  
\end{array}
\end{equation}
This additive resolution was given a moduli-theoretic interpretation in terms of $S^1$-invariant connections on a disk (and Nahm's equations)
in work of Kronheimer \cite{kron-ccad}, Donaldson \cite{Don-bvym} p.114, Kovalev \cite{kovalev-nahm} and Biquard \cite{biquard-nahm} (see especially \cite{biquard-nahm} p.255 for the full resolution picture), and it was shown that $\wt \g$ has a natural holomorphic Poisson structure such that the resolution $\pi$ is the moment map (\cite{biquard-nahm} Th\'eor\`eme 2 (1a)).
This comes down to considering connections of the form $Adz/z$ with $A\in \g$ (and compatible parabolic structures).
The symplectic leaves of $\wt \g$ are the fibres of $\wt \psi$, and are of the form 
$G\times_{B_0} x\lu$ where $x\in \lt=\Lie(T)$ and $\lu$ is the nilradical of $\Lie(B_0)$.

These additive results do not translate into statements for the original (multiplicative) resolution, for example since the exponential map is not surjective in general (e.g. for $\SL_2(\IC)$), and since the centralizer of $\exp(2\pi i x)$ differs from that of $x$ when $x$ reaches the far wall of the Weyl alcove.
Rather, there are natural ``multiplicative'' analogues of the above (additive) symplectic/Poisson statements, as follows.

\begin{thm}[\cite{logahoric}]
For any $t\in T$, the fibre $\wh \cC:=\wt \psi^{-1}(t)\subset \wt G$ is a quasi-Hamiltonian $G$-space with moment map given by the restriction of the resolution map
$$\pi : \wh \cC \to G.$$ 
\end{thm}

In fact a more general statement is proved in \cite{logahoric}
(replacing the Borel subgroup $B_0$ by an arbitrary parabolic $P_0$, and the point $\cC:=\{t\}\subset T$  by an arbitrary conjugacy class $\cC$ of the Levi factor of $P_0$). 
For $G=\GL_n(\IC)$ this was proved earlier by Yamakawa 
\cite{yamakawa-mpa} using quivers.
Our proof proceeds by first constructing a quasi-Hamiltonian $G\times T$-space $\IM=G\times_U B_0$ (\cite{logahoric} Theorem 9) and
then  observing that 
the reduction of $\IM$ by $T$ at the conjugacy class 
$\cC\subset T$ is 
$\wh \cC$. (The spaces $\IM$ are tame analogues of the fission spaces of \S5 above.)
In this approach we can also consider the quotient  $\IM/T$, within the world of quasi-Poisson manifolds \cite{AKM}. This quotient $\IM/T$ is a manifold since the action of $T$ is free, and it is a quasi-Poisson $G$-space for general reasons (from the quasi-Hamiltonian structure on $\IM$), and moreover it is isomorphic to $\wt G$, so we obtain the following:

\begin{cor}
The Grothendieck space $\wt G$ is a quasi-Poisson $G$-space with moment map  $\pi:\wt G\to G.$ 
\end{cor}
 
This is the multiplicative analogue of  the additive Poisson statement above, from \cite{biquard-nahm} Th\'eor\`eme 2. 
The quasi-Poisson bivector is $G$-invariant, has moment map $\pi$, and the leaves of $\wt G$ (in the sense of \cite{AKM} \S9) are the fibres of $\wt \psi$, all analogously to the additive case.
(This contrasts with the Poisson structure on $\wt G$
constructed in \cite{evenslu-groth}.)
Lying behind this is the interpretation of the spaces $\wh \cC$ in terms of monodromy data for logahoric connections (as sketched above, cf. \cite{logahoric} Remark 6, extending Levelt/Simpson for $\GL_n(\IC)$%
), and the realisation that quasi-Hamiltonian/quasi-Poisson geometry {is} the natural geometry of monodromy-type data (in the presence of suitable framings). In brief the element $M\in G$ is the local monodromy, classifying a regular singular meromorphic connection on a $G$-bundle over a punctured disk, and the fibre of $\pi$ over $M$ corresponds to the choice of a logahoric connection extending the regular singular connection across the puncture.

\newpage
\noindent%
{\bf\large9. Dynkin diagrams for isomonodromy systems} \label{sn: slims}

The article \cite{slims} develops a theory of Dynkin diagrams for a large class of isomonodromy systems (i.e. for certain irregular nonabelian Gauss--Manin connections).
In particular we establish a connection between certain moduli spaces of meromorphic connections and a large class of Kac--Moody root systems.
Specifically a class of graphs, the {\em supernova graphs}, is introduced containing
all of the star-shaped graphs as well as all the complete $k$-partite graphs for any integer $k$.
It is shown how one may attach an isomonodromy system to any such graph together with some data on the graph.
Any element of the Weyl group of the Kac--Moody algebra attached to the graph acts on the data and is shown to lift to give an isomorphism between the corresponding isomonodromy systems (often controlling isomonodromic deformations of connections on different rank vector bundles).
Further a characterisation is given (\cite{slims} \S10), in terms of the root system for the Kac--Moody algebra, for exactly when the data determine a non-empty moduli space (this is an additive irregular analogue of the Deligne--Simpson problem).

\noindent{\bf Background.}
In the picture described so far, in previous sections, we may choose a general linear group $G$ and an irregular curve $\Si$ and this determines a \hk manifold as in \cite{wnabh}, which may be viewed in particular as a moduli space $\MDR(\Si)$ of meromorphic connections over $\Si$ with given irregular types.
Then we may vary the irregular curve (in an admissible fashion) over some base $\IB$ to obtain a relative moduli space $\MDR\to \IB$ which has a canonical flat (Ehresmann/nonlinear) connection on it.

Now in certain cases (usually if the underlying algebraic curve is the Riemann sphere) one can do more, and explicitly write down the resulting nonlinear connection as a system of nonlinear differential equations (whence $\IB$ becomes the ``space of times'').
Usually, to get explicit equations, one proceeds by working with a simpler moduli space $\cM^*\subset \MDR(\Si)$
where the underlying vector bundle on the Riemann sphere is holomorphically trivial (this notation goes back to \cite{smid}).
 
If one does this one soon finds however (\cite{Harn94}) that there are different irregular curves (often for different general linear groups)
that lead to the {\em same} system of nonlinear differential equations, cf. \cite{slims} Theorem 1.2.
(Further the underlying moduli spaces of connections are in fact isomorphic---this may be viewed as a failure of the irregular curve analogue of the Torelli problem for these moduli spaces in genus zero.)

This leads to the general question of understanding the isomorphisms (and automorphisms) of such moduli spaces, and the associated isomonodromy systems. 

The simplest examples of such moduli spaces (of complex dimension two) correspond to the six (second-order) Painlev\'e differential equations (these differential equations are explicit expressions of the corresponding nonlinear connections).
In most such cases the moduli spaces $\MDR(\Si)$ are known to coincide with the ``spaces of initial conditions'' constructed explicitly for the Painlev\'e equations by Okamoto \cite{oka-feuil, OkaP24, OkaPVI, OkaPV, Okamoto-dynkin}.
In these works Okamoto also showed that the Painlev\'e equations
admit certain affine Weyl groups of symmetries.
For example the symmetry group of the sixth Painlev\'e equation 
is the affine Weyl group of type $D_4$ (and if one adds diagram automorphisms this extends to affine $F_4$ as used in the work above, in 
\S3, on algebraic solutions).
Similarly the symmetry groups of the fourth and fifth Painlev\'e 
equations are the affine Weyl groups of type $A_2,A_3$ respectively (the others are not simply-laced so will be ignored here for simplicity).

\begin{figure}[h]
	\centering
	\input{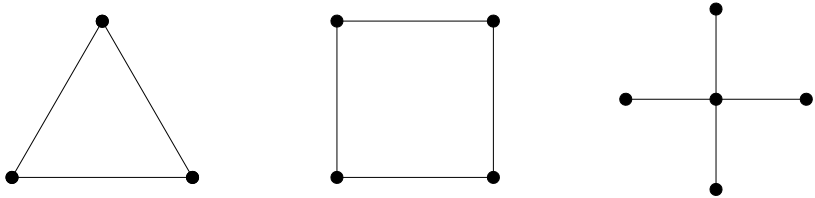}
	\caption{Affine Dynkin diagrams for Painlev\'e equations IV, V and VI.}\label{fig: p456}
\end{figure}

A basic question is thus to understand and extend this link between Dynkin diagrams and the Painlev\'e equations, or isomonodromic deformations more generally. (Needless to say these  affine Weyl groups are not transparent from the moduli problem---for the Painlev\'e equations one starts with certain meromorphic connections on rank two vector bundles on the Riemann sphere.)

Combining the two questions above leads to the following question: can we associate isomonodromy systems to a certain class of graphs, i.e. develop a theory of Dynkin diagrams for isomonodromy systems generalising the above three examples, so that the Weyl group attached to the graph lifts to give automorphisms/isomorphisms?

For example, can we see what is special about the above  three simply-laced affine Dynkin diagrams, that they and no others have associated Painlev\'e equations?

\noindent{\bf Main Results.}  The first step (\cite{quad} Exercise 3) was to notice that the relation Okamoto found between affine Dynkin graphs and
Painlev\'e equations may be understood in a different way.  
In brief Nakajima's theory of {\em quiver varieties} \cite{nakaj-duke94}
gives a way to attach an algebraic variety to a graph and some data on the graph (in fact those of complex dimension two, relevant to Painlev\'e equations, go back at least to 
Kronheimer \cite{Kron.ale}). The observation of \cite{quad} was that for the above Painlev\'e equations the moduli space $\cM^*$ is isomorphic to the quiver variety attached to the corresponding graph. 

Next in \cite{rsode, slims} a higher dimensional version of this observation was established: there are many moduli spaces $\cM^*$ (of meromorphic connections on the trivial bundle on the Riemann sphere) isomorphic to Nakajima quiver varieties, and in fact the full moduli space $\MDR$ is determined by data on the graph.
This extends a relation used by Crawley--Boevey \cite{CB-additiveDS} between star-shaped graphs and Fuchsian systems (the simple pole case).
Whereas quiver varieties are defined for any graph, only for special graphs are there associated moduli spaces of connections. 
The class of simply-laced graphs for which this result holds, 
the {\em supernova graphs}, contains
all of  the complete $k$-partite graphs for any integer $k$ (see \cite{slims}).
\begin{figure}[h]
	\centering
	\input{partite3.pstex_t}
	\caption{Complete $k$-partite graphs from partitions of $N\le 6$}\label{fig: graph table}
\center{\ (omitting the stars $\Ga(1,n)$ and the totally disconnected graphs $\Ga(n)$)}\label{fig: partite-intro}
\end{figure}

Thus, for example, there is no second order Painlev\'e equation attached to the pentagon (the affine $A_4$ Dynkin graph), since it is not a complete $k$-partite graph for any $k$, whereas the square and the triangle are (as is the four pointed star).
Indeed the complete $k$-partite graphs are determined by partitions with $k$ parts and the first few such graphs are as in Figure \ref{fig: graph table}. %

Further, in \cite{slims}, the corresponding isomonodromy systems were written down (attached to any such graph), and given a Hamiltonian formulation.
Most of these systems are new (e.g. they are not included in the work of Jimbo--Miwa et al \cite{JMMS, JMU81}).
It 
was further shown in \cite{slims} how the Weyl group symmetries (of the Kac--Moody algebra with Cartan matrix determined by the graph) lift to relate the isomonodromy systems attached to the spaces of connections. 
This comes about quite naturally by giving a new interpretation of some of the moduli spaces $\cM^*$ (and thus certain quiver varieties) as moduli spaces of presentations of the first Weyl algebra.

Thus, in summary, in certain cases there is an alternative point of view, starting with a graph rather than an irregular curve. 

As a simple application of this way of thinking, by considering hyperbolic Kac--Moody Dynkin graphs (the next simplest class after the affine case), 
\cite{slims} \S11.4 shows there is a family of 
isomonodromy systems (of order $2n$ for any $n\ge 1$) lying over  each of the six Painlev\'e equations.
These are completely different to the well-known ``Painlev\'e hierarchies'' 
and conjecturally related to Hilbert schemes of points on the original two dimensional Painlev\'e moduli spaces.
More precisely (changing complex structure in the \hk family) 
the conjecture of \cite{slims} 
is that the Hilbert scheme of $n$ points on any meromorphic Higgs bundle moduli of complex dimension two, is again a meromorphic Higgs bundle moduli space, and the relation to graphs established in \cite{slims} predicts exactly which higher dimensional moduli space to look at.
(The expected list of such two-dimensional moduli spaces is given in \cite{ihptalk} \S3.2;  the tame cases of this conjecture have apparently been proved recently in \cite{groechenig.hilbert.schemes}).

\

\ 

\ 

\noindent
{\bf\large Other results.}

The articles \cite{roks, quad} solve other long standing problems, perhaps of a more limited interest
(\cite{roks} gives the first conceptual derivation of the famous Regge symmetry of the classical $6j$-symbols, and
\cite{quad} describes the first Lax pairs for the nonlinear additive difference Painlev\'e equations attached to the $E_7$ and $E_8$ root systems \footnote{An $E_6$ Lax pair was obtained earlier in \cite{ArBo-duke}
using the setup of \cite{borodin-it}; \cite{quad} shows that in fact all the additive difference Painlev\'e equations  arise as contiguity relations of systems of linear {\em differential} equations, and so strictly speaking the theory of \cite{ArBo-duke, borodin-it} is unnecessary to understand them.
}).

\renewcommand{\baselinestretch}{1}              %
\normalsize
\bibliographystyle{amsplain}    \label{biby}
\bibliography{../thesis/syr}

\ 

\hfill www.math.ens.fr/$\sim$boalch

\hfill boalch@dma.ens.fr

\end{document}